\documentclass[3p,11pt]{elsarticle}
\usepackage{xcolor}
\usepackage{mathtools}
\usepackage{hyperref} 
\usepackage{setspace}
\usepackage{amsfonts}
\usepackage{amssymb}
\usepackage{amsthm}
\usepackage{bm}
\usepackage{ulem}
\usepackage{caption}
\usepackage{subcaption}
\usepackage{multirow}
\usepackage{array}
\usepackage{comment}
\usepackage{algorithm}
\usepackage{algpseudocode}
\usepackage{soul}
\usepackage{pgfplots}
\pgfplotsset{compat=1.18}
\usepgfplotslibrary{statistics}
\usepackage{booktabs}
\usepackage{todonotes}

\newcommand{\bp}{{\bf p}}
\newcommand{\bn}{{\bf n}}
\newcommand{\ba}{{\bf a}}
\newcommand{\bx}{{\bf x}}
\newcommand{\bg}{{\bf g}}
\newcommand{\bd}{{\bf d}}
\newcommand{\by}{{\bf y}}
\newcommand{\balpha}{{\boldsymbol{\alpha}}}

\newcommand{\stencil}{{\mathcal S}}
\newcommand{\innerpointset}{{\mathcal M_h}}
\newcommand{\ghostpointset}{{\mathcal G_h}}
\newcommand{\allpointset}{{\mathcal A_h}}
\newcommand{\sgn}{{\textrm{sign}}}

\biboptions{numbers,compress}

\newfont{\notapolice}{cmss11}
\newtheoremstyle{myRemark}
{5pt}
{3pt}
{\notapolice}
{}
{\it}
{:}
{.5em}
{}

\newfont{\defpolice}{cmssi11}
\newtheoremstyle{myDefinition}
{5pt}
{3pt}
{\defpolice}
{}
{\bf}
{:}
{.5em}
{}

\theoremstyle{plain}

\theoremstyle{myDefinition}
\newtheorem{defn}{Definition}

\theoremstyle{myRemark}
\newtheorem{remark}{Remark}

\journal{}
\begin{document}

\begin{frontmatter}

\title{Efficient and well-conditioned ghost-point discretization of boundary operators on unfitted domains}

\author[unict]{A. Coco \corref{cor}}
\ead{armando.coco@unict.it}

\author[dei]{A. Coclite}
\ead{alessandro.coclite@poliba.it}

\author[coimbra]{S. L. Clain}
\ead{clain@mat.uc.pt}

\author[minho]{R. M. S. Pereira}
\ead{rmp@math.uminho.pt}

\cortext[cor]{Corresponding author}

\address[unict]{Dipartimento di Matematica, Università degli Studi di Catania, Catania, Italy}

\address[dei]{Dipartimento di Ingegneria Elettrica e dell'Informazione, Politecnico di Bari, Bari, Italy}

\address[coimbra]{Centre of Mathematics (CMUC), University of Coimbra, 3000- 143 Coimbra, Portugal}

\address[minho]{Department of Mathematics, University of Minho, Guimarães, Portugal}

\begin{abstract}

Unfitted boundary methods are widely used to numerically solve partial differential equations (PDEs) on irregular domains, avoiding the computational burden of generating boundary-conforming grids. In the finite-difference framework, structured Cartesian grids offer advantages such as ease of implementation and efficient parallelization, while geometry is represented implicitly, for instance, through level-set functions. In this setting, ghost point methods are commonly employed to enforce boundary conditions by introducing additional relations between interior and ghost nodes. However, constructing these relations becomes challenging for high-order accurate discretizations, which often rely on wide stencils that can reduce computational efficiency and degrade performance in large-scale parallel simulations.

In this work, we investigate alternative ghost-point discretizations based on compact stencils. We introduce a formulation based on a boundary operator that locally approximates the boundary condition near each ghost node, replacing it with linear relations involving both interior and ghost points. The operator is constructed via least-squares reconstruction, allowing flexible stencil configurations while preserving the desired order of accuracy. Several strategies for selecting and adapting compact stencils are proposed, guided by conditioning criteria and iterative refinement procedures to improve global stability. Numerical experiments on various geometries and convection-diffusion regimes demonstrate the effectiveness of the proposed approach, showing that it maintains high accuracy even in the presence of boundary layers and improves stencil compactness and conditioning of the resulting linear systems.
\end{abstract}
\end{frontmatter}

\section{Introduction}\label{sect:intro}
Finite differences on a Cartesian grid offer several advantages, owing to its simple structure, which enables high accuracy while supporting efficient storage and parallel computing.
For problems posed in irregular domains, the computational grid does not conform to the geometry of the domain (unfitted boundary methods).
Instead, the domain is implicitly represented, for instance, by means of a level set function \cite{Osher-Fedkiw:level_set, Gibou:guidelinesPoisson}. Such a method enables the use of Cartesian meshes while allowing the treatment of complex geometries. Moreover, they are particularly efficient in problems involving moving boundaries, since the underlying grid does not need to be recomputed when the geometry evolves in time.

Nevertheless, the unfitted boundary method requires an additional treatment to preserve the global accuracy and stability. To this end, embedded method with ghost points is a popular technique that substitutes the boundary condition on the physical domain with an equivalent thick boundary condition over several layers of external points surrounding the domain~\cite{Peskin:IBM, coclite20191, coclite20204, LeVequeLi:IIM, gibou2002second, Fedkiw:GFM, Gibou:fourth_order, Gibou:fluid_solid}.

On the one hand, the discretization of PDEs yields relations that combine internal points, {\it i.e.}, the nodes that belong to the physical domain, and ghost points. On the other hand, additional relations that combine ghost and internal points have to be found out to close the global system (whose unknowns are the values of the internal and ghost points).
A common first-order approach sets the ghost value equal to that at a near boundary point. To attain higher-order accuracy, the ghost node must be related to interior nodes in a manner that the boundary condition is enforced at a chosen boundary point.
The construction of these relations is the critical issue in a ghost point framework to provide a non-singular system with a small condition number, and several parameters control the design of the additional relations. The principal ingredient to design the complementary relations is the construction of the local stencil, involving ghost and internal nodes together with boundary conditions prescribed at some  point of the physical boundary.

The choice of this boundary point is fundamental and distinguishes the various ghost-point techniques.
Several second-order ghost-point techniques have been proposed in the literature.
In \cite{gibou2002second}, the boundary condition associated with a ghost point is imposed at the intersection between the boundary and the coordinate axes. By contrast, in \cite{coco2013finite, coco2018second, coco2020multigrid} it is imposed at the closest point on the boundary, and this technique is employed in several applications that model complex-shaped boundaries, such as topography in computational volcanology~\cite{coco2016hydro}, monument conservation~\cite{coco2020level}, and moving bubbles in sorption kinetics~\cite{astuto2023finite}.  More recently, a nodal ghost-point formulation on Cartesian meshes has been introduced within the finite element framework~\cite{astuto2024high, astuto2025nodal} and compared with the finite difference scheme in~\cite{Astuto2024Comparison}. This approach provides improved theoretical properties, such as symmetry and coercivity of the discrete operator. However, its extension to high-order accuracy is not straightforward.

Very high-order approximations (third-order or higher)  on unfitted domains require more complex reconstructions within the finite-difference setting. In~\cite{coco2024high}, the authors extend the method of~\cite{coco2013finite} to fourth-order accuracy. We also mention the extension of the extrapolation/interpolation techniques following the normal direction or the axis direction  proposed in ~\cite{Gibou:fourth_order,appelo2012fourth,Baeza2016}. Sixth-order methods are presented in~\cite{fernandez2020very, clain2024very}, where the ghost values are extrapolated solely from nearby interior points and remain decoupled from other ghost points.

In many ghost-point formulations, the stencils adopted for very high-order discretizations are quite extensive and the diameter is larger than 8 to 10 times the characteristic length $\Delta x$, rising several issues: a dramatic discrepancy of the accuracy and the disability to build such a stencil in some situations (corner or wedge). This occurs either because the construction of the discrete operator near the boundary relies on stencils obtained from a square neighbourhood of grid points~\cite{coco2013finite, coco2024high}, or because the stencil does not include other ghost points~\cite{fernandez2020very, clain2024very} (i.e., the ghost equations are decoupled from one another), which leads to the inclusion of interior values that are often far from the boundary. This concern also becomes a limitation in large-scale simulations from a computational point of view, particularly in parallel computing environments, where compact stencils are desirable to reduce communication costs and improve computational efficiency. Moreover, wide stencils may fail to capture boundary effects accurately, such as boundary layers, and this is why compact stencils are desirable to better catch steep gradient solutions.\\

In this paper, we investigate alternative strategies for constructing ghost-point discretizations using more compact stencils. 
We propose a new formulation based on a so-called boundary operator that provides a local discretization of the boundary condition in the vicinity of any ghost point. Such a linear operator mimics the boundary condition up to a given order of accuracy, and substitutes it on the frontier with linear relations over internal and ghost points. In particular, we analyse various configurations of reduced stencils and examine their effects on the resulting discrete operators.  

To this end, we introduce a representation of the boundary condition of any function $\psi$ as a linear combination of the $\psi$ value over the points of a stencil. We say that a reconstruction is $o$-th order if the representation is exact for polynomials up to degree $o-1$. We ensures that the number of points in the stencil is large enough to guarantee that the desired order of accuracy, possibly leading to an underdetermined linear system but with the maximal rank property. We then characterize the quality of the local representation with the condition number of an intermediate non-singular matrix used in the Moore-Penrose system. An additional criterion for assessing the quality of all the representations is the condition number of the global system. Two major criteria control the quality of the reconstruction: the relative position of the points to the target ghost point that constitute the stencil; the choice of the physical point on the boundary where we impose the condition. We then propose several strategies for constructing compact stencils such that we get high quality representations.

As such, the main contributions of this work are: i) the boundary operator formulation that replaces the boundary condition with linear relations between interior and ghost nodes; ii) a least-squares reconstruction framework allowing flexible stencil configurations; iii) conditioning-based stencil selection algorithms ensuring both local and global stability; iv) a cone-based stencil construction strategy that yields compact and well-conditioned stencils.\\

The paper is organized as follows. In Section~\ref{sect:model}, we introduce the general discretization of the PDE and of the boundary conditions, and describe the least-squares approach used to construct the equations for the ghost nodes and to assemble the resulting matrix. In Section~\ref{sect:choice}, we present different stencil strategies for the discretization of the boundary conditions and perform a numerical study to assess their performance and their impact on the conditioning of the resulting linear system. Section~\ref{sect:numexp} is devoted to numerical experiments on different geometries and for various combinations of convection and diffusion coefficients, including tests that highlight boundary layer effects. Finally, conclusions are drawn in Section~\ref{sect:concl}.

%
%
\section{Model and discretization}\label{sect:model}
Let $\Omega\subset \mathbb{R}^2$  be a non-polygonal bounded domain. We seek the solution $\phi:\Omega \rightarrow \mathbb{R}$ of an elliptic problem defined in $\Omega$, complemented with Robin boundary conditions on the 1D smooth closed manifold $\partial\Omega\subset\mathbb{R}^2$.

The problem then reads
\begin{equation}
\begin{cases}
-k\Delta \phi({\bf x}) + U({\bf x})\cdot \nabla \phi({\bf x}) = f({\bf x}) \qquad &\text{for} \quad {\bf x} \in\Omega \\
a_D({\bf x})\phi({\bf x}) + a_N({\bf x})\nabla \phi({\bf x})\cdot {\bf n}({\bf x}) = g({\bf x}) \qquad &\text{for} \quad {\bf x} \in\partial\Omega
\end{cases}\, ,
\label{generic}
\end{equation}
with $k>0$ the diffusion coefficient; $U({\bf x})\in \mathbb{R}^2$ the assigned velocity, and $f({\bf x})$ the source term. Moreover, $a_D({\bf x}),\, a_N({\bf x})$ stand for the Dirichlet and Neumann contribution parameters with ${\bf n({\bf x})}$ the outward normal vector at any point ${\bf x}\in\partial \Omega$, while $g({\bf x})$ is a regular function.
%
%

\subsection{Notations}
Given a spacing parameter $h$, we define the Cartesian grid characterized by point $(x_i,y_j)$ with $x_i = ih$ and $y_j = jh$.  
The set $\innerpointset$ gathers the points $\bx_k$ of the grid that belong to $\Omega$, with indices $k=1,\cdots,N_I=|\innerpointset|$. The set $\ghostpointset$ corresponds to the ghost points of the grid $\bx_k\notin\Omega$ constituted of one or two layers of points around $\Omega$ in function of the expected accuracy (2nd, or 4th, respectively). Layer 1 is constituted of the horizontally and vertically adjacent points to $\innerpointset$, and we repeat the procedure for the second and third layers. 

We complete the identification of the ghost points with 
index $k=N_I+1\cdots, N_I+N_G$, with $N_G=|\ghostpointset|$. Set $\allpointset=\innerpointset\cup\ghostpointset$ is the set of {\it active points}, that are the points where we aim at defining an approximation $\phi_k\approx \phi(x_k)$. We store the values in $\Phi=(\Phi_I,\Phi_G)$, with $\Phi_I$ and $\Phi_G$ the vectors associated to $\innerpointset$ and $\ghostpointset$ respectively. We define a mapping $k\to (i,j)=\big (i(k),j(k)\big)$ such that $\bx_k=(x_i,y_j)$, $k=1,\cdots,N_I+N_G$.

\subsection{Inner operator discretization }
Let us consider an $o$-th order finite-difference discretization of the differential operators in \eqref{generic}, then for any $\bx_k=(x_i,y_j) \in\innerpointset$ the discrete problem reads
\begin{equation}
A^o_x[\Phi] + A^o_y[\Phi] = f_I\, ,
\label{discrete}
\end{equation}
where $A^o_x[\Phi]$ and $A^o_y[\Phi]$ represent the discrete left-hand side operator, while $f_I\in\mathbb{R}^{N_I}$ is the vector collecting the values $f(\bx_k)$.
Specifically, we design $N_I$ algebraic equations for $\bx_k\in\innerpointset$
\begin{equation}
A^o_x[\Phi]_{k} + A^o_y[\Phi]_{k} = f_{k}=f(\bx_k)\,,
\end{equation}
where $A^o_x[\Phi]_h$ and $A^o_y[\Phi]_h$ stand for the discrete left-hand side operator while $f_I\in\mathbb{R}^{N_I}$ the vector collects the values $f(\bx_k)$.
In this manuscript, we focus on a fourth-order centred scheme for the Laplace operator as well as for the gradient; as such, $A^4_x[\Phi]_{k}$ is obtained as
\begin{equation}
\begin{split}
A^4_x[\Phi]_{k} = &-k\frac{-\frac{1}{12}\phi_{i-2,j}+\frac{4}{3}\phi_{i-1,j}-\frac{5}{2}\phi_{i,j}+\frac{4}{3}\phi_{i+1,j}-\frac{1}{12}\phi_{i+2,j}}{h^2} \\ 
&+ u\frac{\frac{1}{12}\phi_{i-2,j}-\frac{2}{3}\phi_{i-1,j}+\frac{2}{3}\phi_{i+1,j}-\frac{1}{12}\phi_{i+2,j}}{2h}
\end{split}\, ,
\label{FDx}
\end{equation}
with $\bx_k=(x_i,y_j)$ and analogously $A^4_y[\Phi]_{k}$ reads
\begin{equation}
\begin{split}
A^4_y[\Phi]_{k} = &-k\frac{-\frac{1}{12}\phi_{i,j-2}+\frac{4}{3}\phi_{i,j-1}-\frac{5}{2}\phi_{i,j}+\frac{4}{3}\phi_{i,j+1}-\frac{1}{12}\phi_{i,j+2}}{h^2} \\ 
&+ u\frac{\frac{1}{12}\phi_{i,j-2}-\frac{2}{3}\phi_{i,j-1}+\frac{2}{3}\phi_{i,j+1}-\frac{1}{12}\phi_{i,j+2}}{2h}.
\end{split}\,
\label{FDy}
\end{equation}
\begin{remark}
While point $\bx_k=(x_i,y_j)\in\innerpointset$, this is not necessarily the case for the other points $\bx_\ell$ that are involved in the formula. Nevertheless, the ghost points' set $\ghostpointset$ has to be designed to guarantee that $\bx_\ell\in\allpointset$.
\end{remark}

Consequently, Eq~\eqref{discrete} corresponds to an $N_I\times N$ linear system ($N_I$ algebraic equations with $N$ unknowns) involving all points of $\innerpointset$ as well as all points of $\ghostpointset$, formulated under the matrix form:  
\begin{equation}
A_{II}\Phi_I+A_{IG}\Phi_G=f_I\, ,
\label{discrete_MI}
\end{equation}
where $A_{II}$ and $A_{IG}$ label the discrete operators acting on interior or ghost points. 
%
%

\subsection{Discrete boundary operator}
To close the system \eqref{discrete_MI}, one has to provide $N_G$ additional algebraic equations corresponding to the entries of the vector $\Phi_G$. To this end, let $\bx_k$ be a ghost point and $\bp_k\in\partial \Omega$ a point close to $\bx_k$. For any smooth function $\psi$, we have
$$
\mathfrak{B}[\psi]({\bp_k}) = a_D({\bp_k})\psi({\bp_k})+a_N({\bp_k})\nabla \psi({\bp_k})\cdot \bn({\bp_k}).
$$
Now, let  $\stencil^k$ be a stencil of size $N_k=|\stencil^k|$ with $\bx_k\in \stencil^k$ and constituted of ghost points or inner points close to $\bx_k$. The procedure to build efficient and well-conditioned $\stencil^k$ will be the topic of the next section. 

At that stage, for any regular function $\psi$, we introduce the functional
\begin{equation}
\ba \in\mathbb R^{N_k}\to \mathcal{E}^k(\ba;\stencil^k)[\psi]=\left ( \sum_{\ell\in\stencil^k} a_\ell \psi(\bx_\ell)\right )-\mathfrak B[\psi](\bp_k)\, ,
\label{stencil:func}
\end{equation}
where $\ba$ collects the coefficients $a_\ell$ associated with the points of the stencil $\stencil^k$. On the other hand, $\mathbb P^{o-1}$ is the polynomial space of degree lower or equal to\ $o-1$, equipped with the basis of $\displaystyle N_o=\frac{o(o+1)}{2}$ independent polynomials 
$$
\psi(\bx;\bx_k,\balpha)=(\bx-\bx_k)^\balpha,\ |\balpha|\leq o-1,
$$
with $\balpha=(\alpha_x,\alpha_y)$, $|\balpha|=\alpha_x+\alpha_y$ and $(\bx-\bx_k)^\balpha=(x-x_k)^{\alpha_x}(y-y_k)^{\alpha_y}$.
In this paper, we aim to achieve fourth-order accuracy for both the solution and its gradient, under Dirichlet as well as Neumann boundary conditions. While fourth-order accuracy of the solution, along with Dirichlet boundary conditions, would require only $o = 4$, to attain the desired accuracy for the gradient and in the presence of mixed boundary conditions, we use $o = 5$ (then $\mathbb P^4$ polynomial space), corresponding to $N_o=15$.

To determine the vector $\ba$, we state that the functional has to satisfy the following constraints:
$$
\mathcal{E}^k(\ba;\stencil^k)[\psi(\bx;\bx_k,\balpha)]=0,\quad |\balpha|\leq o-1.
$$
We then have $N_k$ unknown coefficients with $N_o$ constraints and reads $B^k\ba=\bg^k$ with
\begin{equation}\label{def_matrix_B}
B^k_{\ell,\balpha}=(\bx_\ell-\bx_k)^\balpha,\qquad \bg^k_\balpha=\mathfrak B[(\bx-\bx_k)^\balpha](\bp_k),\quad \ell\in\mathcal S^k,\ |\balpha|\leq o-1.
\end{equation}
\begin{remark}
At the implementation level, we need an indexation $\balpha\to m(\balpha)$ with $m$ ranging from $1$ to $N_o$ to manipulate the matrix and carry out the computations.
\end{remark}
$\blacktriangleright$ If $N_k<N_o$, we do not have enough degrees of freedom to guarantee the desired order of accuracy, which disqualifies such a situation. One has to choose a stencil large enough such that $N_k\geq N_o$.

$\blacktriangleright$ If $N_k=N_o$ and assuming that the square matrix $B^k$ is non-singular, we then obtain a unique vector $\ba^k$ to represent the boundary condition operator at point $\bp_k$.

$\blacktriangleright$ If $N_k>N_o$, we have an underdetermined system regarding the coefficients, and we use the least squares method under constraints. More precisely, we seek a vector $\ba^k$ such that:    
$$
\ba^k=\arg \min_{\ba \in\mathbb R^{N_k}} \frac{1}{2} \ba^t \ba
$$
under the linear constraints $\mathcal{E}^k(\ba;\stencil^k)[\psi(\bx;\bx_k,\balpha)]=0$, $|\balpha|\leq o-1$. 
\begin{remark}\label{remark:pk}
A key note is to observe that $\ba^k=\ba^k(\bp_k)$ since the determination of $\ba^k$ as a minimizer of \eqref{stencil:func} depends on the choice of point $\bp_k\in\partial \Omega$. Thus, we have an additional degree of freedom from the choice of $\bp_k$ that will be exploited to provide the "best"  boundary operator.
\end{remark}

As a consequence, given $\bp_k$ and the stencil $\mathcal S_k$, we determine the vector $\ba^k$, and the discrete boundary operator at the ghost point $\bx_k\in\ghostpointset$ reads:
\begin{equation}\label{eq:akl}
\mathfrak B_k(\Phi)=\sum_{\bx_\ell\in\mathcal S_k} a^k_\ell \phi_\ell.
\end{equation}
At last, the new equation associated to the ghost point $\bx_k$ gives
\begin{equation}\label{eq:problpha}
\mathfrak B_k(\Phi)=g({\bp}_k),\quad  k=N_I+1,\dots,N_I+N_G\,.
\end{equation}
In that way, we define the $N_G$ equations from the relations \eqref{eq:problpha}, that we can rewrite under the compact form:

\begin{equation}\label{eq:Ng_eqs}
A_{GI} \Phi_I + A_{GG} \Phi_{G}=F_G\, .
\end{equation}
with $F_G$ the vector of entries $g({\bp}_k)$.

Assembling the two systems Eqs~\eqref{discrete_MI} and~\eqref{eq:Ng_eqs} provides the global discrete $N\times N$ problem 
\begin{equation}\label{eq:systems}
\left [ \begin{array}{cc}
A_{II} & A_{IG}\\
A_{GI} & A_{GG}
\end{array} \right ]
\left [ \begin{array}{c}
\Phi_I\\ \Phi_G
\end{array} \right ] =
\left [ \begin{array}{c}
F_I\\F_G
\end{array} \right ]\, ,
\end{equation}
where matrices $A_{II}$, $A_{IG}$ come from the interior operator discretization whereas matrices $A_{GI}$ and $A_{GG}$ derives from the boundary operator discretization.
\begin{remark}
The linear system \eqref{eq:Ng_eqs} has to be of maximal rank $N_G\times N$ to provide an eligible global system. In fact, we need to prescribe more constraints about the stencil design and the choice of $\bp_k$ to control the condition number, both locally (matrix $B_k$) and globally (the global linear system).
\end{remark}

\subsection{The level-set function and the boundary projection}\label{sect:levelset}
Let the domain $\Omega$ be represented by a level–set function 
$\varphi(\mathbf{x})$, which is negative inside $\Omega$, positive outside, and vanishes on the boundary:
\[
    \Omega = \{\,\mathbf{x} \in \mathbb{R}^2 : \varphi(\mathbf{x}) < 0 \,\}, 
    \qquad
    \partial \Omega = \{\,\mathbf{x} \in \mathbb{R}^2 :  \varphi(\mathbf{x}) = 0 \,\}.
\]
The outward unit normal to the boundary is given by
$\mathbf{n} = \nabla \varphi / \left| \nabla \varphi \right| $,
where $\varphi$ is assumed to be explicitly available. Given a ghost point $\mathbf{x}_k$, its orthogonal projection onto the boundary is computed as an iterative procedure as given in \cite{coco2024high}.

\section{The choice of the stencils}\label{sect:choice}
\subsection{Stencil quality and condition number}
Stability and accuracy of the approximation derive from two key stages. On a local scale, the stencil selection results in the conditioning of the local matrix, which connects the ghost point with the boundary data and interior points. At a large scale, the global matrix, which encompasses all the interior and ghost points, is generated from the local matrices. The global condition number is essential to obtain accurate approximations for general boundary configurations. Consequently, it is imperative to address the two issues for any non-polygonal frontier and boundary conditions and develop a robust and efficient algorithm to select stencils.
\vskip 1em

We first introduce the local criterion to assess the quality of the stencil and the associated discrete boundary operator, characterized by the vector $\ba^k$.
\begin{defn}\label{defn:admissible}
Let $\stencil^k$ be a stencil associated to a ghost point $\bx^k$.  The stencil is {\bf admissible} if the matrix $B^k$ given by \eqref{def_matrix_B} satisfies the maximum rank property.
\end{defn}
If we satisfy the maximum rank property then $B^{k\dagger}=(B^k)^\top B^k$ is non-singular and one has $\ba^k=B^{k\ddag}\bg^k$, with $B^{\ddag}=(B^{k\dagger})^{-1}(B^k)^\top$ the Moore-Penrose  matrix.
\begin{defn}\label{defn:localwellcond}
An admissible stencil is said to be {\bf locally well conditioned} if $B^{k\dagger}$ is a well-conditioned matrix.
Accordingly, we adopt the criterion
$\chi(B^{k\dagger}) < \Lambda_{\mathrm{LOC}}$,
where $\chi(B^{k\dagger})$ denotes the condition number and $\Lambda_{\mathrm{LOC}}$ is a prescribed tolerance.
\end{defn}

On the other hand, we deal with the global matrix using a second criterion to control stability and accuracy. A heuristic argument concerns the condition number of the square matrix $A_{GG}$ that connects the ghost points between themselves. We expect to have a well-conditioned matrix $A_{GG}$ if each diagonal entry is not too small compared to the off-diagonal entries in the same row (in modulus). This requirement leads to the following definition.
\begin{defn}\label{defn:wellcond}
The influence of the stencil $\stencil^k$ on the global conditioning of the matrix $A_{GG}$ is measured by the value 
\begin{equation}\label{eqn:ratio}
R_k = \frac{1}{\vert a_k^k\vert} \max_{\bx_\ell \in \mathcal G_k}^{\ell \neq k} \vert a^k_\ell\vert,
\end{equation}
where $a^k_\ell$ are given in~\eqref{eq:akl}.
An admissible stencil $\stencil^k$ is said {\bf globally well-conditioned} with respect to the ghost point $\bx_k$ and to the boundary point $\bp_k$ if one has,
\[
R_k < \Lambda_{\mathrm{GLO}},
\]
where $\Lambda_{\mathrm{GLO}}$ is a prescribed tolerance.
\end{defn}
\begin{remark}
We recall from Remark~\ref{remark:pk} that $\ba^k=\ba^k(\bp_k)$ for any $\bx_k \in\ghostpointset$, where $\bp_k\in \partial \Omega$ is close to the point $\bx_\ell$. Therefore, the quantity $R_k$ deeply relies on the choice of the points where we prescribe the boundary condition on the physical interface.
In summary, the admissibility of a stencil only depends on the relative positions of the stencil grid points, whereas the well-conditioning criterion depends on the mutual position of stencil grid points and the boundary point $\bp_k$ where the conditions are prescribed.
\end{remark}

\begin{remark}
While the non-singularity of the global matrix is a mandatory requirement to provide a solution, the construction of local and global well-conditioned stencils are desirable properties that reduce computational effort and provide better accuracy. 
\end{remark}

\vskip 1em

In conclusion, the main challenge is the design of algorithms that generate optimal stencil in such a way that we satisfy the two criteria given by Definitions \ref{defn:localwellcond} and \ref{defn:wellcond}.
To this end, we first investigate triangular stencils, which automatically satisfy Definition~\ref{defn:localwellcond} by construction. Within this class, we identify a configuration that also fulfills Definition~\ref{defn:wellcond}. However, this stencil may exhibit undesirable properties, such as a large diameter, which can lead to poor spectral resolution. Consequently, an alternative stencil-generation procedure is required.
To address this limitation, we propose a strategy that does not rely on triangular stencils. Specifically, we first construct a buffer of candidate points and iteratively select points until Definition~\ref{defn:localwellcond} is satisfied. We then verify whether Definition~\ref{defn:wellcond} holds. If this condition is not met, we identify a critical point responsible for the issue and we replace it with a new candidate from the buffer. This process is repeated until all required conditions are satisfied. As we will demonstrate, the resulting stencil is significantly more compact.

\begin{remark}\label{rmk:expl}
Note that, if $\bx_k$ is the unique ghost-point belonging to $\stencil^k$, the problem is explicit, giving that the approximation of $\phi_G$ at $\bx_k$ is obtained through an affine boundary operator: 
\begin{equation}
(\bx_k\in\stencil^k, \textrm{boundary condition on }\bp_k) \to \phi_G(\bx_k).
\end{equation}
Such a procedure is a highly parallelizable computational framework.
\end{remark}

\subsection{The stencil construction}\label{sect:ST}
As anticipate at the end of the previous section, we explore different constructions to provide a suitable stencil that will be analyzed in terms of the following requirements:\\
\noindent [1:] admissible stencil;\\
\noindent [2:] locally well conditioned stencil;\\
\noindent [3:] globally well conditioned stencil;\\
\noindent [4:] compact stencil.\\[0.5em]
To formalize the different strategies, we introduce additional notation.
Let $\bx_k = (x_k,y_k)$ be a ghost-point and $\bp_k=(p_x^k,p_y^k)$ the associated point on the physical boundary. Vector $\bd_k=\bx_k-\bp_k=(d_x^k,d_y^k)$ is the displacement vector between $\bx_k$ and $\bp_k$, and $p>1$ an integer that will control the stencil size $N_k=|\stencil^k|$ according to the target accuracy we intend to reach.

\subsubsection{Stencil based on triangle shape}
\begin{figure}[ht]
\begin{center}
    \begin{subfigure}{0.32\textwidth}
        \centering
        \includegraphics[width=\textwidth]{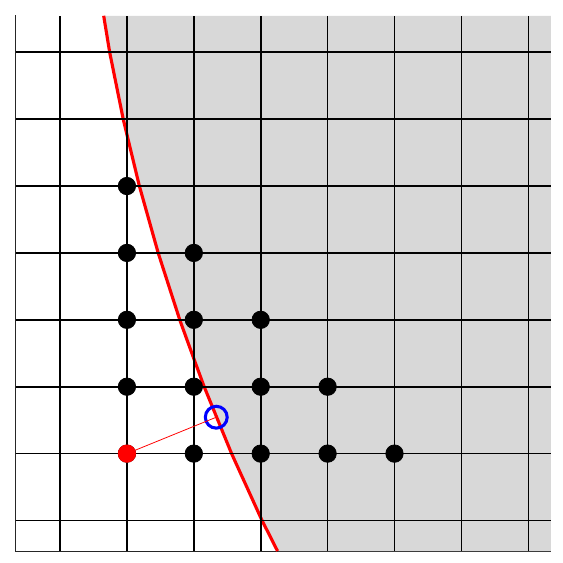}
        \caption{}
    \end{subfigure}
    \begin{subfigure}{0.32\textwidth}
        \centering
        \includegraphics[width=\textwidth]{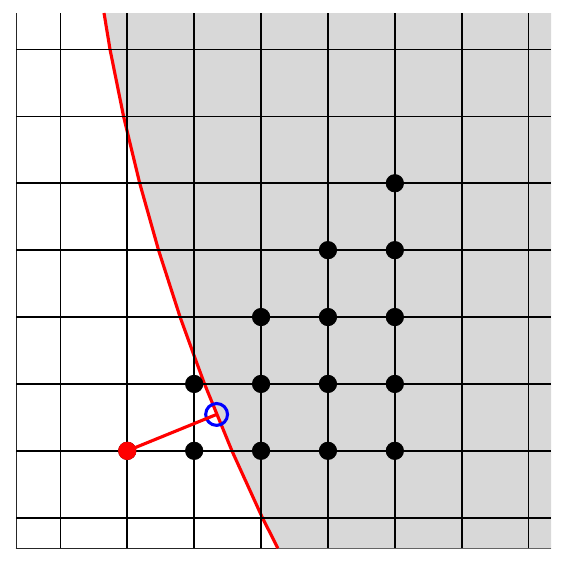}
        \caption{}
    \end{subfigure}
    \begin{subfigure}{0.32\textwidth}
        \centering
        \includegraphics[width=\textwidth]{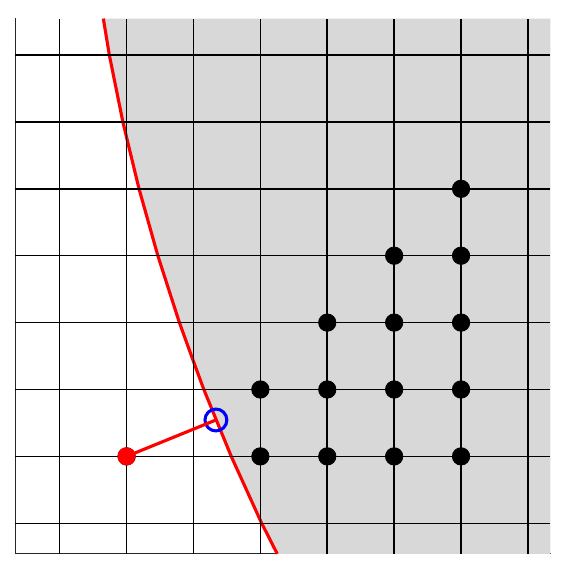}
        \caption{}
    \end{subfigure}
\end{center}
\caption{Stencils S1 (a), S2 (b) and S3 (c). The red point denotes the ghost node ${\bf x}_k$ from which the stencil is constructed, while the blue circle indicates the corresponding boundary point ${\bf p}_k$.}
\label{fig:S123}
\end{figure}

We first consider the stencil design based on points belonging to a triangular zone of the grid. We aim to produce a compact stencil that balances the number of horizontal and vertical displacements to ensure homogeneity in both directions. 
\paragraph{Stencil S1: External Right-Triangular Stencil}
It is a right triangle-shaped stencil with the ghost point positioned at the vertex of the right angle (see Figure~\ref{fig:S123}.a):
\begin{equation}
\label{st:trinojump}
\stencil^k = \left\{ (x_k + l\,\sgn(d_x^k),y_k + m\,\sgn(d_y^k)), \; l,m=0,\ldots,p, \; l+m \leq p  \right\}\, .
\end{equation}
Such a right-angle triangular stencil satisfies requirements 1 and 2. 
\paragraph{Stencil S2: Internal Right-Triangular Stencil}
To reduce the number of ghost points (other than $\bx_k$) in this right-triangular stencil, the vertex of the right angle is chosen as the innermost internal point $(x_R^k,y_R^k)$ (see Figure~\ref{fig:S123}.b):
\begin{equation}
\label{eq:modprob}
\left\{
\begin{aligned}
(x_R^k,y_R^k)&=(x_k+p\,\sgn(d_x^k),y_k), \quad \text{if } |x_k-p_x^k| \geq |y_k-p_y^k|\, , \\
(x_R^k,y_R^k)&=(x_k,y_k+p\,\sgn(d_y^k)), \quad \text{otherwise.}
\end{aligned}
\right.
\end{equation}
Then:
\begin{equation}\label{st:trifrominsidenojump1}
\stencil^k = \left\{ (x_R^k - l\,\sgn(d_x^k),y_k + m\,\sgn(d_y^k)), \; l,m=0,\ldots,p, \; l+m \leq p  \right\}\, ,
\end{equation}
if $|x_k-p_x^k| \geq |y_k-p_y^k|$, or, otherwise,
\begin{equation}
\label{st:trifrominsidenojump2}
\stencil^k = \left\{ (x_k + l\,\sgn(d_x^k),y_k - m\,\sgn(d_y^k)), \; l,m=0,\ldots,p, \; l+m \leq p  \right\}.
\end{equation}
This stencil satisfies requirements 1 and 2. However, if $|\bx_k-\bp_k|>h$, the magnitude of the coefficients associated with the ghost-point would result in a much smaller value than those associated with other points of the stencil. Therefore, requirement 3 is not satisfied.
\paragraph{Stencil S3: Internal Right-Triangular Stencil with Jump}
To satisfy requirement 3 when considering an internal right-triangular stencil S2, additional internal grid points have to be included rather than ghost points.
At this scope, if $|\bx_k-\bp_k|>h$, the diameter of the stencil is enlarged by one column and a line of ghost points is substituted with a line on internal points (see Figures~\ref{fig:S123}.c).
Notably, this stencil contains no ghost points other than ${\bx}_k$; therefore, the problem would result as explicit (see Remark~\ref{rmk:expl}), and requirements 1, 2, and 3 are satisfied.

\subsubsection{Stencil construction based on a cone (Stencil S4)}
\begin{figure}[ht]
\begin{center}
    \begin{subfigure}{0.32\textwidth}
        \centering
        \includegraphics[width=\textwidth]{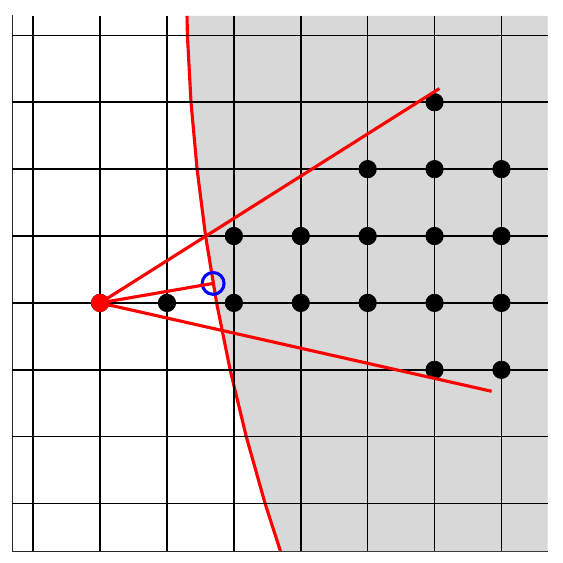}
        \caption{}
    \end{subfigure}
    \begin{subfigure}{0.32\textwidth}
        \centering
        \includegraphics[width=\textwidth]{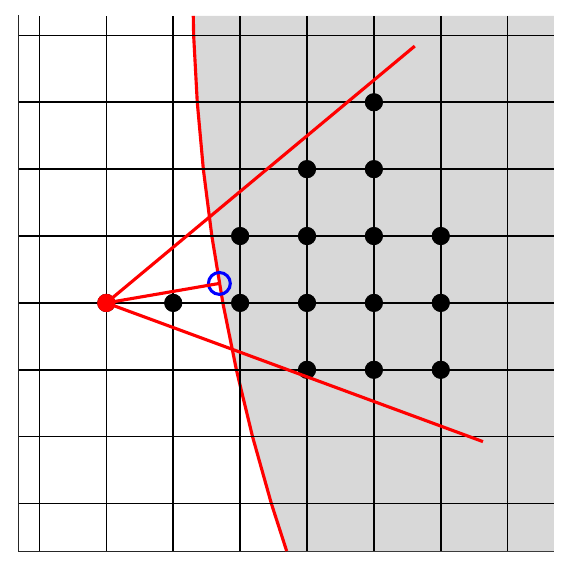}
        \caption{}
    \end{subfigure}
    \begin{subfigure}{0.32\textwidth}
        \centering
        \includegraphics[width=\textwidth]{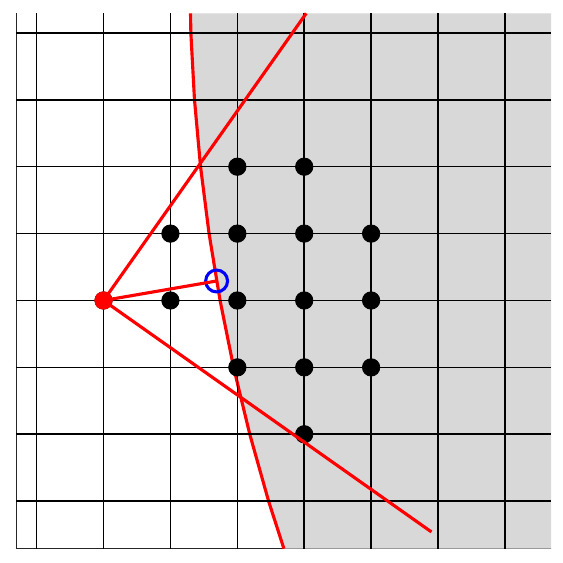}
        \caption{}
    \end{subfigure}
    \begin{subfigure}{0.32\textwidth}
        \centering
        \includegraphics[width=\textwidth]{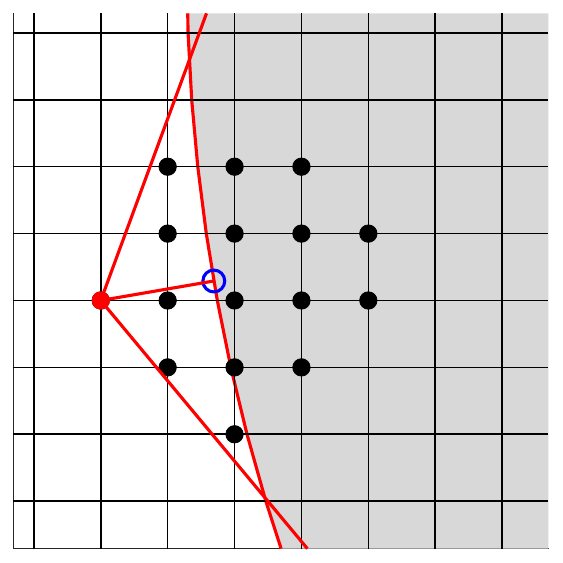}
        \caption{}
    \end{subfigure}
    \begin{subfigure}{0.32\textwidth}
        \centering
        \includegraphics[width=\textwidth]{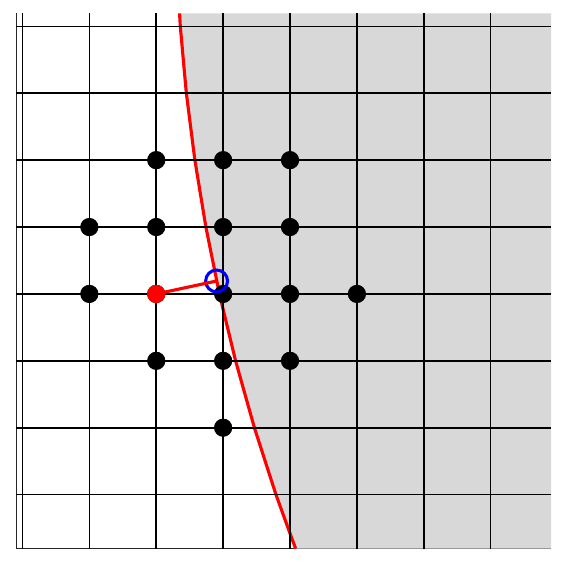}
        \caption{}
    \end{subfigure}
\end{center}
\caption{Stencils S4.1, with $\vartheta = 45^\circ$ (a), $60^\circ$ (b), $90^\circ$ (c), $120^\circ$ (d), $360^\circ$ (e). The red point denotes the ghost node ${\bf x}_k$ from which the stencil is constructed, while the blue circle indicates the corresponding boundary point ${\bf p}_k$.}
\label{fig:S2}
\end{figure}

\begin{figure}[ht]
\begin{center}
\includegraphics[width=0.45\textwidth]{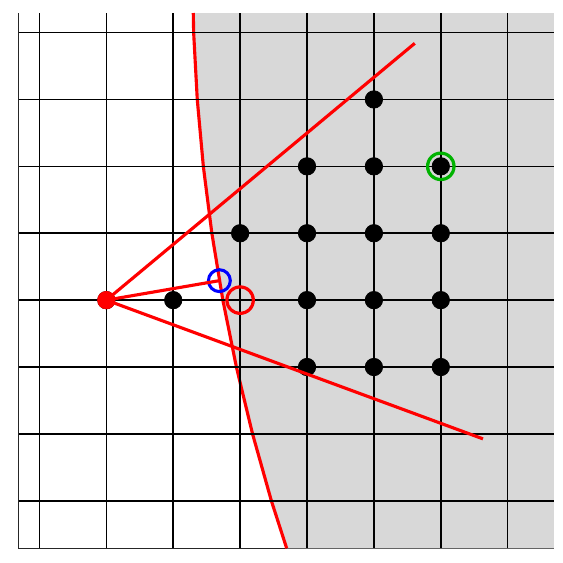}
\includegraphics[width=0.45\textwidth]{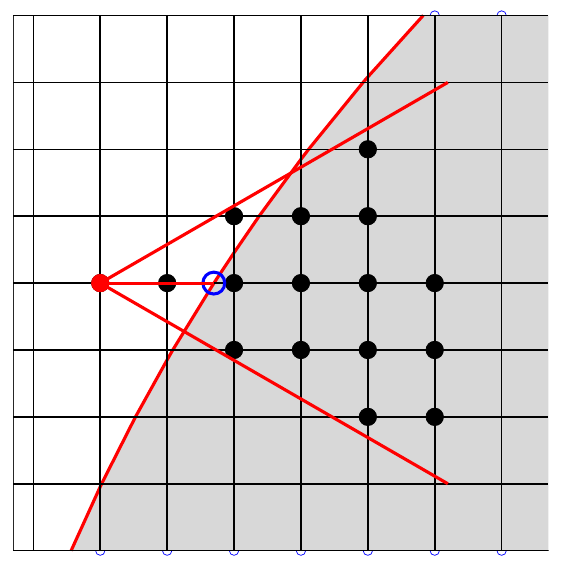}
\end{center}
\caption{
Stencils S4.2 (left) and S4.3 (right) with $\vartheta = 60^\circ$. The red point denotes the ghost node ${\bf x}_k$ from which the stencil is constructed, while the blue circle indicates the corresponding boundary point ${\bf p}_k$. For the Stencil S4.2 (left), the red circle denotes the node removed from the stencil, which is replaced by the node indicated by the green circle.}
\label{fig:stencils_S42andS43}
\end{figure}

We consider the half-cone $\mathcal C(\bx_k,\bp_k,\vartheta)$ centred at the point $\bx_k$ of aperture $\vartheta$ with semi-axis given by the vector $\bx_k\bp_k$. Then we constitute the stencil by picking-up $N_o$ active nodes (ghost and interior nodes) $\by_k\in \mathcal C(\bx_k,\bp_k,\vartheta)$, close enough to $\bp_k$, {\it i.e.} the point that belongs to the physical boundary. We display in Figure \ref{fig:S2} several situations for characteristic angles of the cone.
Unlike the S1–S3 stencils, this stencil is not necessarily admissible according to Definition~\ref{defn:admissible}.
To address such a concern, we employ three successive strategies, applied hierarchically, to generate three stencils: S4.1, S4.2 and S4.3.
\begin{description}
\item[S4.1:] Incrementally add candidate points to the stencil until it becomes both admissible and locally well-conditioned, according to Def.~\ref{defn:localwellcond} using a tolerance of $\Lambda_\text{LOC} = 10^6$.
\item[S4.2:] Compute the stencil S4.1 first. If the stencil is not globally well-conditioned according to Def.~\ref{defn:wellcond} and with a tolerance of $\Lambda_\text{GLO} = 10$, remove the point associated with the largest coefficient $a^k_\ell$ (see Definition~\ref{defn:wellcond}), and replace it with the nearest available candidate point (see Figure~\ref{fig:stencils_S42andS43}, left panel). This procedure is repeated until the stencil becomes globally well-conditioned, up to a maximum of three iterations.
\item[S4.3:] Compute the stencil S4.2 first. If the stencil remains not globally well-conditioned according to Def.~\ref{defn:wellcond} and with a tolerance of $\Lambda_\text{GLO} = 10$, modify the collar point $\bp_k$ by taking the horizontal or vertical projection from the ghost point onto the boundary (choosing the closer one if both projections exist). Then, apply the strategy described above to construct the S4.2 stencil using this modified boundary point (see Figure \ref{fig:stencils_S42andS43}, right panel).
\end{description}

The strategies are summarized with the Algorithm~\ref{alg:S4} in \ref{sec::appA}.

\section{Numerical Experiments}\label{sect:numexp}
This section presents a series of benchmarks designed to evaluate the advantages and limitations of the stencil design strategies introduced above. The domain $\Omega$ is immersed in a square $[-1,1]\times[-1,1]\subset \mathbb{R}^2$ equipped with a Cartesian grid of characteristic length $h$. The set $\mathcal{M}_h$ gathers the interior nodes, whereas $\mathcal{G}_h$ is composed by the ghost nodes of the grid. 
To determine the order of convergence, we consider 13 levels of grid refinement with characteristic mesh size $h=\frac{2}{N}$, where $N$ increases exponentially. Specifically, $N=$ 160, 176, 194, 213, 234, 258, 283, 312, 343, 377, 415, 456, and 502. Numerical errors are computed using the $L^1$ and $L^\infty$ norm, as detailed in the next section. 

\subsection{Errors and orders computation}
Let us denote the analytical solution $\phi^a(\bx)$ and the numerical approximation  $\phi^h(\bx)$. We quantify the numerical errors in both the solution and its gradient. In particular, we consider the relative $L_1$ error and the $L_\infty$ errors defined on the discrete mesh $\mathcal{M}_h$ as:
\begin{equation}\label{eq:errorSol}
\|e\|_{L_1}
=
\frac{\sum_{\bx_i \in \mathcal{M}_h}
\left|\phi^h(\bx_i)-\phi^a(\bx_i)\right|}
{\sum_{\bx_i \in \mathcal{M}_h}
\left|\phi^a(\bx_i)\right|},
\qquad
\|e\|_{L_\infty}
=
\max_{{\bf x}_i \in \mathcal{M}_h}
\left|\phi^h(\bx_i)-\phi^a(\bx_i)\right|.
\end{equation}
The relative $L_1$ error and the $L_\infty$ error of the gradient are defined as:
\begin{equation}\label{eq:errorGrad}
\|\nabla^h e\|_{L_1}
=
\frac{\sum_{\bx_i \in \mathcal{M}_h}
\left|\nabla^h \phi^h(\bx_i)-\nabla \phi^a(\bx_i)\right|}
{\sum_{\bx_i \in \mathcal{M}_h}
\left|\nabla \phi^a(\bx_i)\right|},
\qquad
\|\nabla^h e\|_{L_\infty}
=
\max_{\bx_i \in \mathcal{M}_h}
\left|\nabla^h\phi^h(\bx_i)-\nabla \phi^a(\bx_i)\right|,
\end{equation}
where the pointwise gradient error magnitude is given by
\[
\left|\nabla^h \phi^h(\bx_i)-\nabla \phi^a(\bx_i)\right|
=
\left[
\left(\partial_x^h \phi^h(\bx_i)-\partial_x \phi^a(\bx_i)\right)^2
+
\left(\partial_y^h \phi^h(\bx_i)-\partial_y \phi^a(\bx_i)\right)^2
\right]^{1/2}.
\]
Here, $\partial_x^h$ and $\partial_y^h$ denote fourth-order finite-difference reconstructions of the gradient computed from the numerical solution $\phi^h$. For example, the discrete $x$-derivative is given by
\[
\partial_x^h \phi^h(x,y)
=
\frac{-\phi^h(x+2h,y)+8\phi^h(x+h,y)-8\phi^h(x-h,y)+\phi^h(x-2h,y)}{12h},
\]
with an analogous expression used for $\partial_y^h \phi^h$.

\subsection{The homogeneous Poisson problem on an annulus domain}\label{sect:test_hom}
We first present a simple benchmark to assess stencil design strategies applied to Dirichlet and/or Neumann conditions. The goal is to select the algorithm that provides the best stencil in terms of accuracy and stability. To this end, we consider an annulus domain $\Omega\subset \mathbb{R}^2$ with radii $R_1=\frac{\sqrt{5}}{5}$ and $R_2=\frac{\sqrt{3}}{2}$, namely $\Omega = \left\{ \bx \in \mathbb{R}^2 \;\middle|\; R_1 < |\bx| < R_2 \right\}$. The inner and outer boundaries are $\partial_i\Omega = \left\{ \bx \in \mathbb{R}^2 \;\middle|\; |\bx| = R_1\right\}$ and $\partial_o\Omega = \left\{ \bx \in \mathbb{R}^2 \;\middle|\; |\bx| = R_2\right\}$, respectively.  The boundary of the domain is then given by $\partial\Omega = \partial_i\Omega \cup \partial_o\Omega$.

For the sake of simplicity, problem~\eqref{generic} is simplified to the Laplace operator with $k = 1$, $U(\bx) = (0,0)$, $f(\bx) = 0$ in $\Omega$. Moreover, for the boundary condition we use $a_D(\bx) = 1$ and $a_N(\bx) = 0$ on $\partial_i\Omega$, and $a_N(\bx)= 1$ and $a_D(\bx) = 0$ on $\partial_o\Omega$. The problem then reads
\begin{equation}
\begin{cases}
-\Delta \phi(\bx) = 0 \qquad &\text{for} \ \ \bx\in\Omega, \\
\phi(\bx) = 0 \qquad &\text{for} \ \ \bx\in\partial_i\Omega,\\
\nabla \phi(\bx)\cdot \bn(\bx) = 1 \qquad &\text{for} \ \ \bx\in\partial_o\Omega.
\end{cases}
\label{poisson}
\end{equation}

Note that the rotational invariance of the parameters and geometry implies that the analytical solution is radial (depending only on $|\bx|$) and reads $\phi^a(\bx)=a+b\log(|\bx|)$. Enforcing the boundary conditions, we have: 
\begin{equation}
a+b\log(R_1)=0, \ \ \frac{b}{R_2}=1\, .
\end{equation}
Then, $b=R_2$ and $a=-R_2\log(R_1)$, leading to 
\begin{equation}
\phi^a({\bf x})=R_2\log\left(\frac{|\bx|}{R_1}\right)\, .
\label{solP}
\end{equation}

\begin{figure}[ht]
\begin{center}
    \begin{subfigure}{0.32\textwidth}
        \centering
        \includegraphics[width=\textwidth]{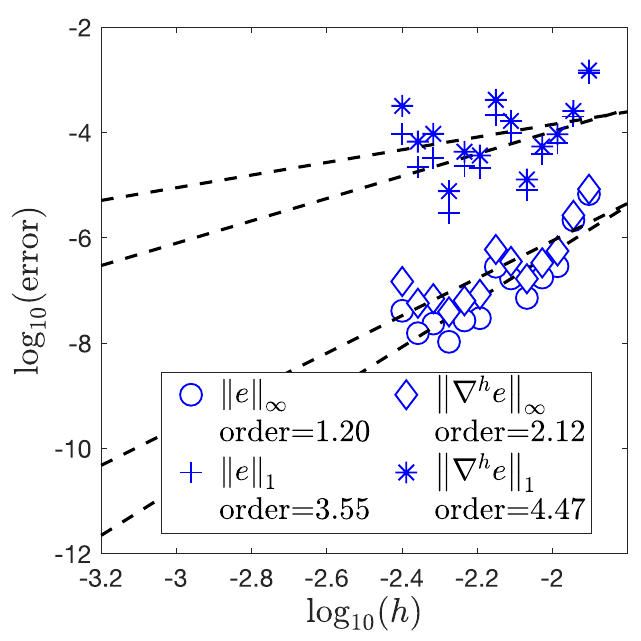}
        \caption{}
    \end{subfigure}
    \begin{subfigure}{0.32\textwidth}
        \centering
        \includegraphics[width=\textwidth]{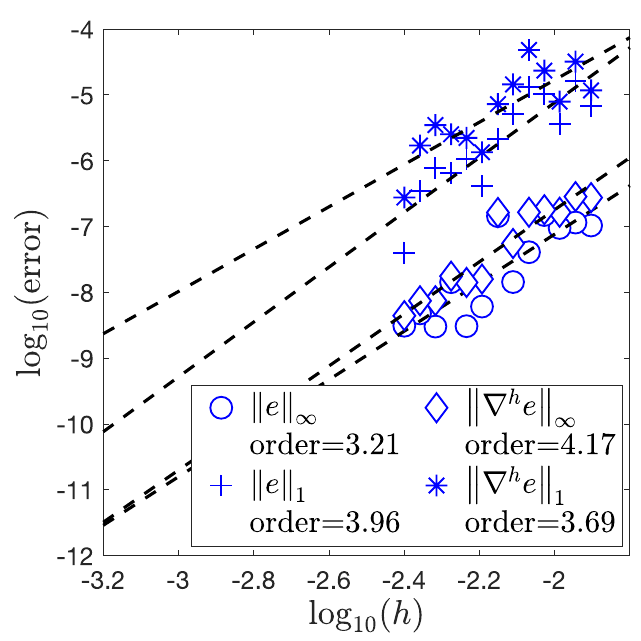}
        \caption{}
    \end{subfigure}
    \begin{subfigure}{0.32\textwidth}
        \centering
        \includegraphics[width=\textwidth]{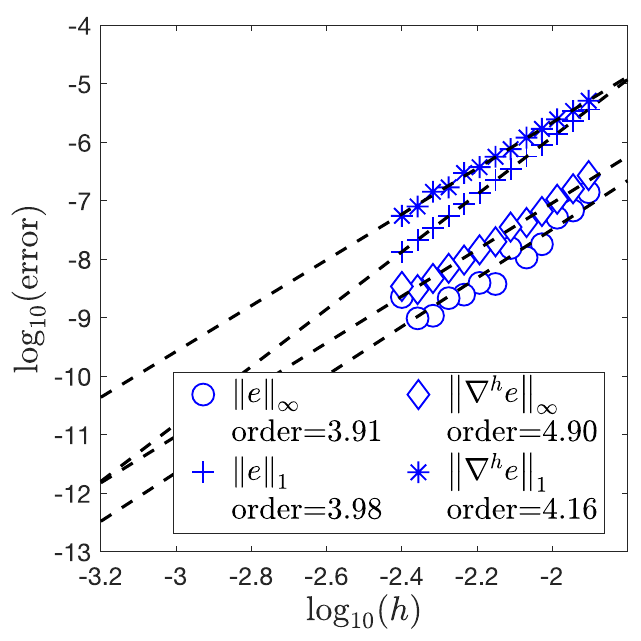}
        \caption{}
    \end{subfigure}
\end{center}
\caption{Numerical errors for the test described in Section~\ref{sect:test_hom}.
We use the stencils S1 (a), S2 (b), and S3 (c).}
\label{fig:BF_S123}
\end{figure}

We first analyze the impact of the S1, S2 and S3 stencil reconstructions for the simple Poisson problem. We display in Fig.~\ref{fig:BF_S123} (left panel) the errors in norm $L^1$ and $L^\infty$ by using the S1 algorithm. Computations have been carried out over the family of grids and we draw in blue the average decay from which we deduce the order method. The middle panel and the right panel correspond to the S2 and S3 algorithms, respectively.
Stencil S1 leads to a scheme with low accuracy, mainly due to the treatment of ghost points located far from the boundary, where the resulting stencil is not guaranteed to be globally well-conditioned.
Stencil S2 improves the error decay, although the convergence is not strictly monotonic.
A more consistent decrease is achieved with stencil S3, which is globally well-conditioned; however, the $L^\infty$ norm of the gradient of the solution does not reach the expected fourth-order accuracy.

\begin{figure}[ht]
\begin{center}
\includegraphics[width=0.3\textwidth]{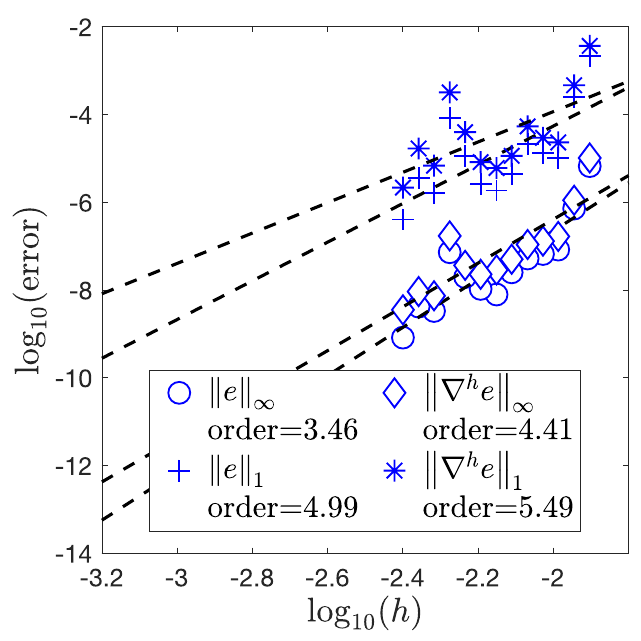}
\includegraphics[width=0.3\textwidth]{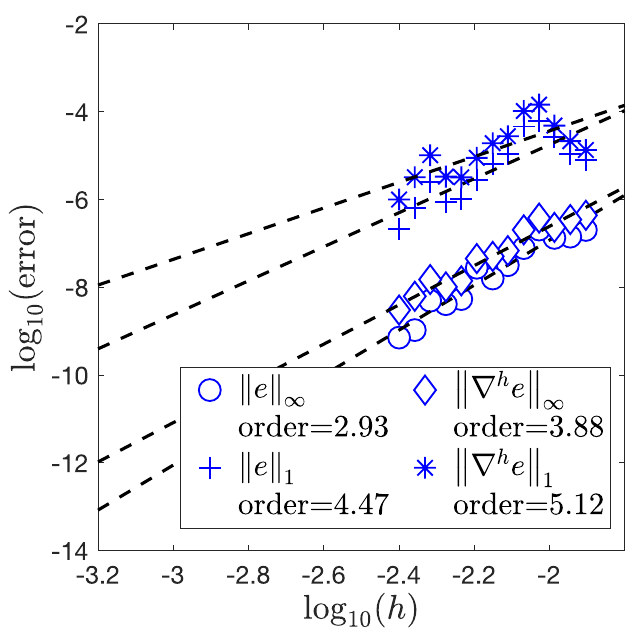}
\includegraphics[width=0.3\textwidth]{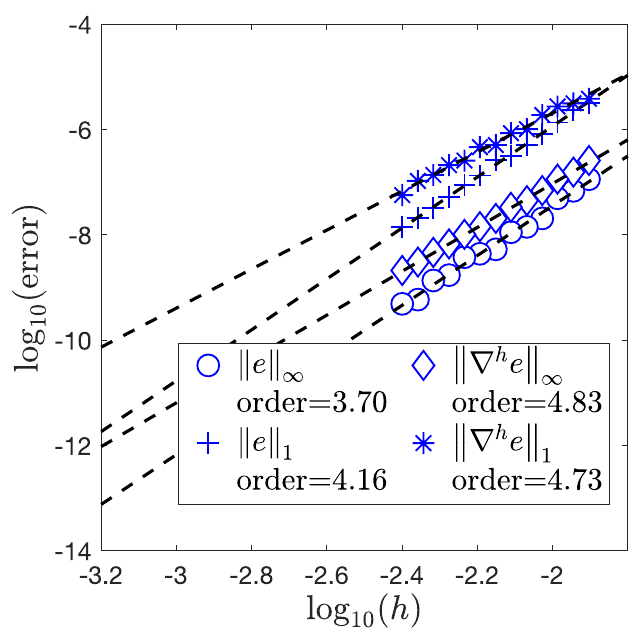}
\end{center}
\caption{
Numerical errors for the test described in Section~\ref{sect:test_hom}.
We use the stencils S4.1 (left), S4.2 (middle), and S4.3 (right), with a vertex angle of $\vartheta = 60^\circ$.
}
\label{fig:STR123}
\end{figure}

We now consider the S4 reconstruction and examine the three variants for a fixed cone angle of $\vartheta = 60^\circ$.
Fig.~\ref{fig:STR123} shows the errors and convergence curves for the first (left panel), second (middle panel), and third (right panel) variants across the different grids.
Despite stencil S4.1 yields reasonable results for the $L^1$ norms of both the solution and its gradient, the $L^\infty$ norm does not decrease monotonically and therefore we discard such a version. Stencil S4.2 shows a noticeable improvement, while the best overall performance is achieved with S4.3.

To provide a more in-depth analysis, we report in Table~\ref{tab:length} the stencil sizes for the different variations. The first column lists the number of points in each stencil, while the remaining columns report the number of stencils for the corresponding sizes. We observe that none of the three stencils (S4.1, S4.2 and S4.3) requires more than 19 points, and the vast majority of stencils contain at most 17 points (recall that $N_o=15$). 
Thus, admissibility is attained without a significant enlargement of the stencil.

\begin{table}[ht]
\centering
\begin{tabular}{c cc cc cc}
\toprule
 & \multicolumn{6}{c}{Stencil length} \\
\cmidrule(lr){2-7}
\# of grid points
& \multicolumn{2}{c}{S4.1}
& \multicolumn{2}{c}{S4.2}
& \multicolumn{2}{c}{S4.3} \\
\cmidrule(lr){2-3} \cmidrule(lr){4-5} \cmidrule(lr){6-7}
 & Count & (\%) & Count & (\%) & Count & (\%) \\
\midrule
15 & 912 & 24.46 & 912 & 24.46 & 808 & 21.67 \\
16 & 1348 & 36.16 & 1188 & 31.87 & 1284 & 34.44 \\
17 & 1032 & 27.68 & 1112 & 29.83 & 1192 & 31.97 \\
18 & 404 & 10.84 & 484 & 12.98 & 412 & 11.05 \\
19 & 32 & 0.86 & 32 & 0.86 & 32 & 0.86 \\
\bottomrule
\end{tabular}
\caption{Stencil size (number of points per stencil) for the three stencils S4.1, S4.2, and S4.3, with $N=502$.}
\label{tab:length}
\end{table}

\begin{figure}
\phantom{m}\hskip -2em
\begin{minipage}{0.45\textwidth}
\begin{tikzpicture}
\begin{axis}[
    boxplot/draw direction=y,
    xtick={1,2,3},
    ylabel={$\log_{10}(\chi(B^\dagger))$},
    xticklabels={S4.1, S4.2, S4.3},
    axis lines=box,
    cycle list={{black}},          
    boxplot/every box/.style={draw=black, fill=gray!20},
    boxplot/every median/.style={red},
    boxplot/every whisker/.style={black, thick},
    boxplot/box extend=0.25, 
]

\addplot+[
    boxplot prepared={
        lower whisker=4.0515,
        lower quartile=4.1227,
        median=4.631,
        upper quartile=4.7636,
        upper whisker=5.0735
    }
] coordinates {};

\addplot+[
    boxplot prepared={
        lower whisker=4.0515,
        lower quartile=4.2247,
        median=4.7097,
        upper quartile=4.7881,
        upper whisker=5.8272
    }
] coordinates {};

\addplot+[
    boxplot prepared={
        lower whisker=4.0515,
        lower quartile=4.2247,
        median=4.6706,
        upper quartile=4.7636,
        upper whisker=5.3298
    }
] coordinates {};

\end{axis}
\end{tikzpicture}
\end{minipage}
\hspace{0.05\textwidth} 
\begin{minipage}{0.45\textwidth}
\begin{tikzpicture}
\begin{axis}[
    boxplot/draw direction=y,
    xtick={1,2,3},
    ylabel={$\log_{10}(R_k)$},
    xticklabels={S4.1, S4.2, S4.3},
    axis lines=box,
    cycle list={{black}},          
    boxplot/every box/.style={draw=black, fill=gray!20},
    boxplot/every median/.style={red},
    boxplot/every whisker/.style={black, thick},
    boxplot/box extend=0.25, 
]

\addplot+[
    boxplot prepared={
        lower whisker=0,
        lower quartile=0.044535,
        median=0.44195,
        upper quartile=0.9843,
        upper whisker=3.3349
    }
] coordinates {};

\addplot+[
    boxplot prepared={
        lower whisker=0,
        lower quartile=0.044535,
        median=0.4378,
        upper quartile=0.80644,
        upper whisker=3.0802
    }
] coordinates {};

\addplot+[
    boxplot prepared={
        lower whisker=0,
        lower quartile=0.044535,
        median=0.4378,
        upper quartile=0.79162,
        upper whisker=1.4987
    }
] coordinates {};

\end{axis}
\end{tikzpicture}
\end{minipage}
\caption{
Boxplots comparing local and global conditioning for the three stencils S4.1, S4.2, and S4.3, with $N=502$.
Left: distribution of condition numbers of the local matrices $B^\dagger$ on a logarithmic scale.
Right: distribution of the value $R_k$ (Eq.~\ref{eqn:ratio}) on a logarithmic scale.
Boxes indicate the interquartile range, the central line denotes the median, and whiskers extend to the minimum and maximum values.
}
\label{fig:boxplots}
\end{figure}
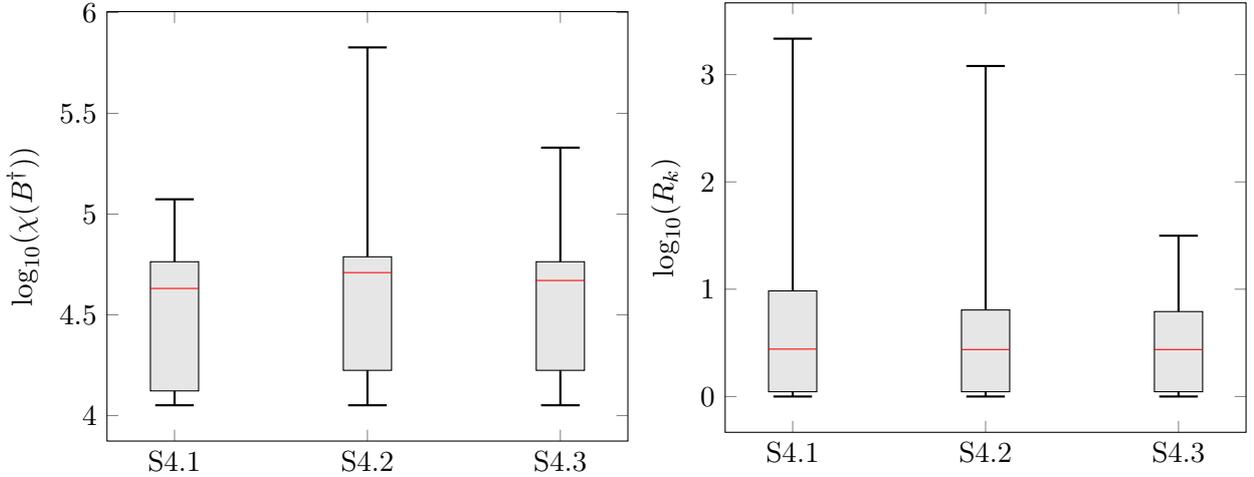

Figure \ref{fig:boxplots} compares the three stencils (S4.1, S4.2, and S4.3) with respect to Requirements 1-3, using $N=502$.
The left panel shows a box plot of the condition numbers on a logarithmic scale of the matrices 
$B_k$ (as defined in Definition \ref{defn:localwellcond}). Stencil S4.1 produces some highly ill-conditioned stencils, with condition numbers exceeding $3 \times 10^5$, whereas stencils S4.2 and S4.3 substantially improve conditioning.
Finally, the right panel reports the values of $R_k$ (see Definition \ref{defn:wellcond}). Stencil S4.3 consistently yields stencils that satisfy the global well-conditioning criterion, performing best among the three stencils in this regard.

\begin{figure}[ht]
\begin{center}
    \begin{subfigure}{0.32\textwidth}
        \centering
        \includegraphics[width=\textwidth]{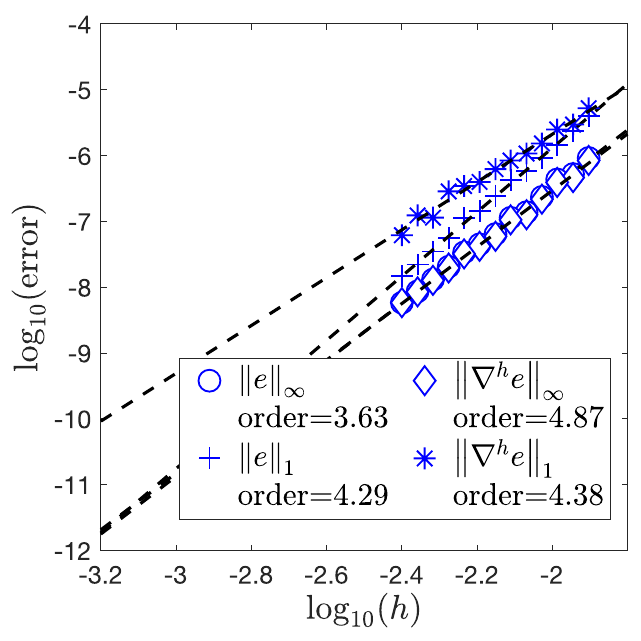}
        \caption{}
    \end{subfigure}
    \begin{subfigure}{0.32\textwidth}
        \centering
        \includegraphics[width=\textwidth]{figures/bestfit_STRATEGY3_CONE60_nstencil15.pdf}
        \caption{}
    \end{subfigure}
    \begin{subfigure}{0.32\textwidth}
        \centering
        \includegraphics[width=\textwidth]{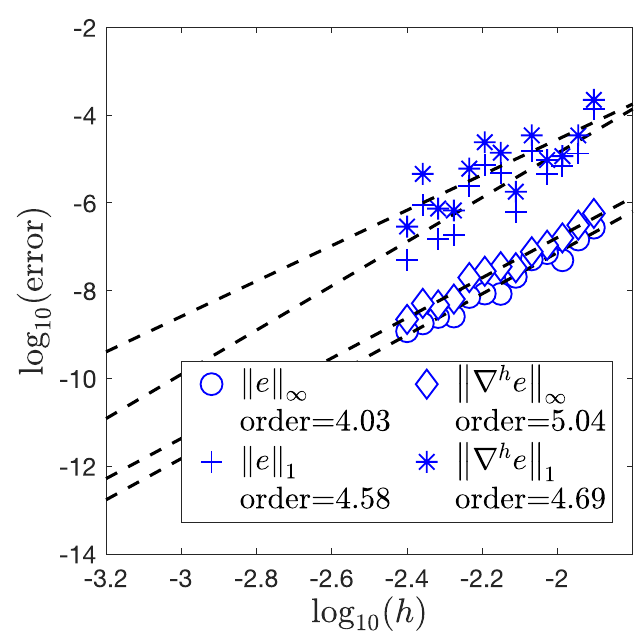}
        \caption{}
    \end{subfigure}
    \begin{subfigure}{0.32\textwidth}
        \centering
        \includegraphics[width=\textwidth]{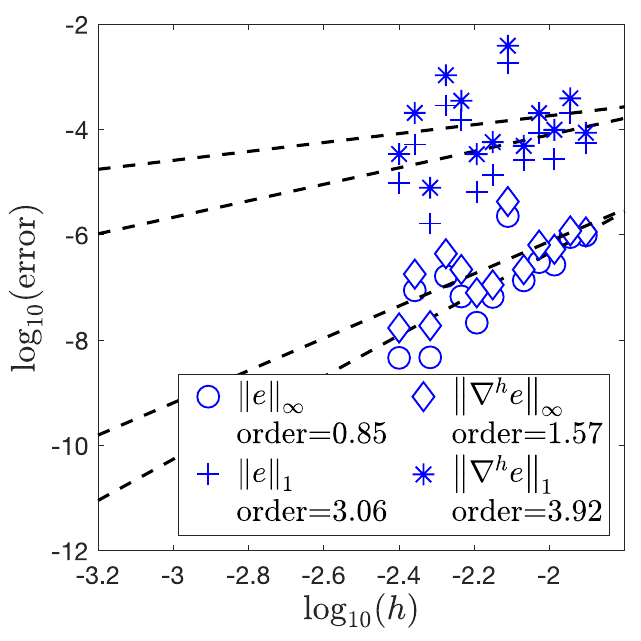}
        \caption{}
    \end{subfigure}
    \begin{subfigure}{0.32\textwidth}
        \centering
        \includegraphics[width=\textwidth]{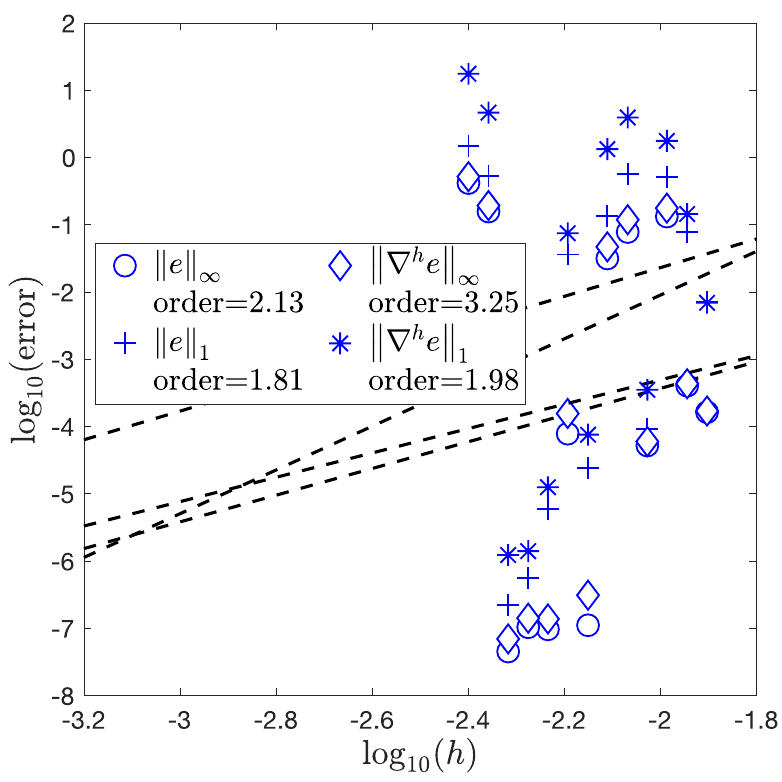}
        \caption{}
    \end{subfigure}
\end{center}
\caption{Numerical errors for the test described in Section~\ref{sect:test_hom}.
We use the stencil S4.3 and a vertex angle of $\vartheta = 45^\circ$ (a), $60^\circ$ (b), $90^\circ$ (c), $120^\circ$ (d), $360^\circ$ (e). See Fig.~\ref{fig:S2} for the respective stencils.
}
\label{fig:BF_S2}
\end{figure}

The choice of the cone angle $\vartheta$ introduces an additional degree of freedom. We therefore analyze the stencil quality as a function of the aperture. To this end, we employ the S4.3 algorithm to construct stencils for different angles, namely  $\vartheta =30,\,45,\,90,\,120,\,360^\circ$, reproducing the error and convergence curves in Fig.~\ref{fig:BF_S2}. We clearly achieve the  best results with $\vartheta = 45^\circ$ and $\vartheta = 60^\circ$.

\vskip 1em 

\begin{figure}
\begin{center}
\begin{tikzpicture}
\begin{axis}[
    boxplot/draw direction=y,
    xtick={1,2,3},
    ylabel={Stencil diameter},
    xticklabels={Approach A, Approach B, Approach C},
    axis lines=box,
    every axis plot/.append style={black},
    cycle list={{black}},          
    boxplot/every box/.style={draw=black, fill=gray!20},
    boxplot/every median/.style={red},
    boxplot/every whisker/.style={black, thick},
    boxplot/box extend=0.25, 
]

\addplot+[
    boxplot prepared={
        lower whisker=5.6569,
        lower quartile=5.6569,
        median=6.4031,
        upper quartile=6.4031,
        upper whisker=7.2111
    }
] coordinates {};

\addplot+[
    boxplot prepared={
        lower whisker=5.831,
        lower quartile=5.831,
        median=6.4031,
        upper quartile=6.7082,
        upper whisker=7.6158
    }
] coordinates {};

\addplot+[
    boxplot prepared={
        lower whisker=5,
        lower quartile=5.099,
        median=5.3852,
        upper quartile=5.6569,
        upper whisker=6.0828
    }
] coordinates {};

\end{axis}
\end{tikzpicture}
\end{center}
\caption{
Boxplots of stencil diameters (maximum pairwise distance between stencil points divided by the spatial step $h$) for the three most accurate approaches: S3 (Approach A), and S4.3
for $\vartheta = 45^\circ$ (Approach B) and $\vartheta = 60^\circ$ (Approach C). We choose $N=502$.
Boxes indicate the interquartile range, the central line denotes the median, and whiskers extend to the minimum and maximum values.
}
\label{fig:histo_diam}
\end{figure}
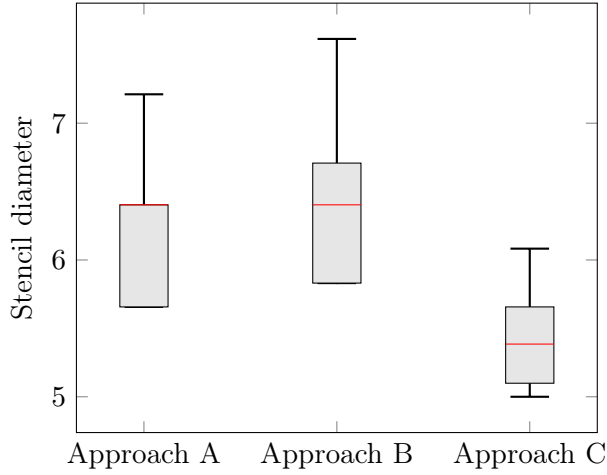

As an intermediate conclusion, we observe that the best results are obtained with the algorithms S3 (Approach A), S4.2 with $\vartheta = 45^\circ$ (Approach B) and S4.3 with $\vartheta = 60^\circ$ (Approach C). To compare the three approaches, we now assess Requirement 4. Figure \ref{fig:histo_diam} provides estimates of the stencil normalized diameters, defined as the stencil diameter divided by the step size $h$. Approach C clearly provides the most compact stencils, and therefore corresponds to the best approach for Requirement 4. In conclusion, from this point onwards, we adopt Approach C (stencil S4.3 with a vertex angle of $\vartheta = 60^\circ$) for all subsequent tests.

\subsection{Poisson problem and complex domains}\label{sect:otherDomains}
\begin{figure}[ht]
\begin{center}
\includegraphics[width=0.30\textwidth]{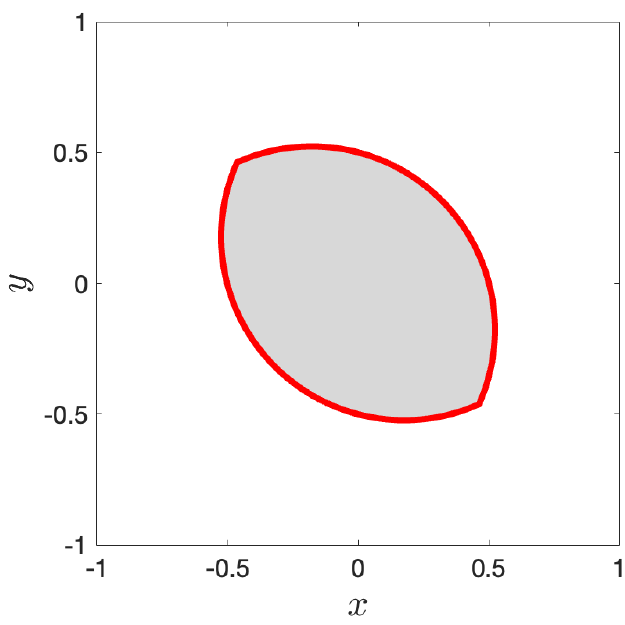}
\includegraphics[width=0.30\textwidth]{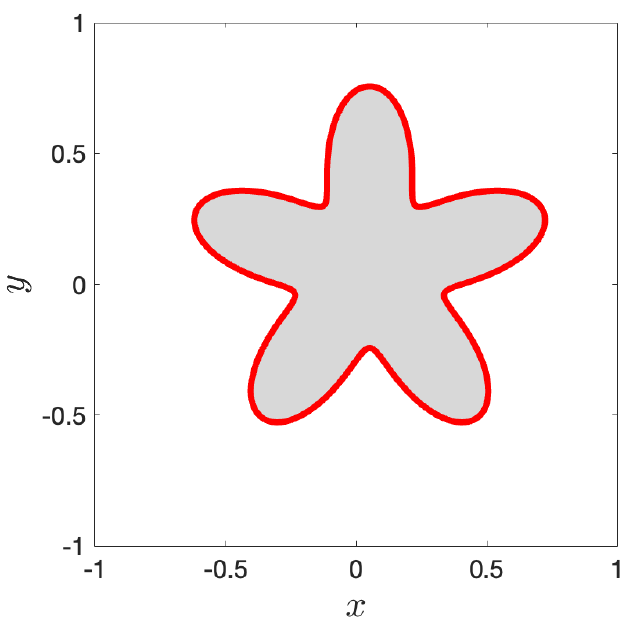}
\includegraphics[width=0.30\textwidth]{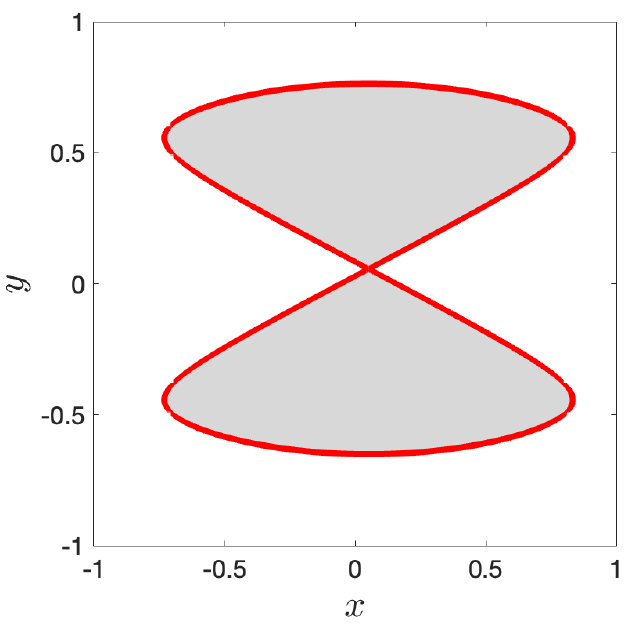}
\end{center}
\caption{Complex geometries that do not fit with the Cartesian grid: the rotated leaf (left), the flower (middle), and the hourglass (right).}
	\label{fig:shape_domains}
\end{figure}

We address the ability of the numerical method to handle complex domains that do not fit with the Cartesian grids. In particular, we consider domains with different types of curvature (locally concave or convex), or presenting some cusps or corners. This leads to challenges in prescribing the correct boundary conditions at the discrete level, while preserving the optimal order. We assess the accuracy of the proposed method across different geometries (taken from~\cite{astuto2025nodal}): the rotated leaf, the flower, and the hourglass. These three domains are shown in Fig.~\ref{fig:shape_domains} and represented by a level set function, as described in Section~\ref{sect:levelset}.

$\bullet$ The rotated leaf. We use the following level-set function:
\begin{align*}
&\widetilde x_1 = -0.25,\, x_1 = \widetilde x_1\cos(\pi/4),\,
y_1 = \widetilde x_1\sin(\pi/4),
\quad R_1 = \sqrt{(x-x_1)^2+(y-y_1)^2} \\
&\widetilde x_2 = 0.25,\, x_2 = \widetilde x_2\cos(\pi/4),\,
y_2 = \widetilde x_2\sin(\pi/4),
\quad R_2 = \sqrt{(x-x_2)^2+(y-y_2)^2} \\       
&r_0 = 0.7,\quad  \varphi_1 = R_1-r_0, \quad \varphi_2 = R_2-r_0, \quad \varphi = \max\{\varphi_1,\varphi_2\}.
\end{align*}

$\bullet$ The flower-shaped domain. The level-set function is:
\begin{align*}
&X = x-0.03\sqrt 3, \quad  Y = y-0.04\sqrt 2, \quad R = \sqrt{X^2+Y^2} \\
&\varphi = R - 0.5 - \frac{Y^5 +5X^4Y-10X^2Y^3}{5R^5}.   
\end{align*}

$\bullet$ The hourglass-shaped domain. It is a domain with a saddle point, whose level-set function is:
\begin{align*}
    X = x-0.03\sqrt{3}, \quad Y = y-0.04\sqrt{2}, \quad \varphi = 256\,Y^4-16\,X^4-128\,Y^2 + 36\,X^2.
\end{align*}

We solve the Poisson equation with source term and mixed boundary conditions:
\begin{equation}
\begin{cases}
-\Delta \phi({\bf x}) = f \qquad &\text{for} \ \ {\bf x}\in\Omega \\
\phi({\bf x}) = g_D \qquad &\text{for} \ \ {\bf x}\in\Gamma_D\\
\nabla \phi({\bf x})\cdot {\bf n}({\bf x}) = g_N \qquad &\text{for} \ \ {\bf x}\in\Gamma_N
\end{cases}\,
\label{poisson}
\end{equation}
where $\Gamma_D = \left\{ (x,y) \in \partial \Omega \text{ such that } x\geq 0 \right\}$ and 
$\Gamma_N = \partial \Omega \backslash \Gamma_D$.
We choose $f$, $g_D$ and $g_N$ in such a way that the analytical solution is $\phi^a(x,y) = \sin(2x)\sin(5y)$, namely 
$f = - \Delta \phi^a$, $g_D = \phi^a$ and $g_N = \nabla \phi^a \cdot \frac{\nabla \varphi}{|\nabla \varphi|}$.
\vskip 1em

\begin{figure}[ht]
\begin{center}
\includegraphics[width=0.3\textwidth]{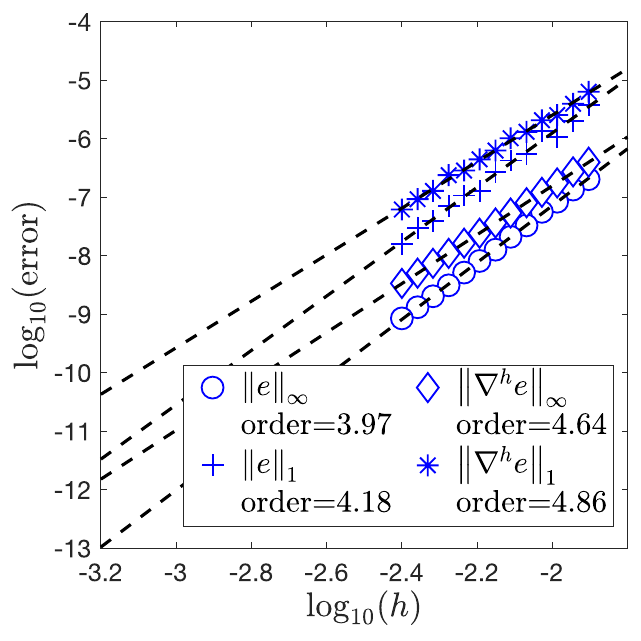}
\includegraphics[width=0.3\textwidth]{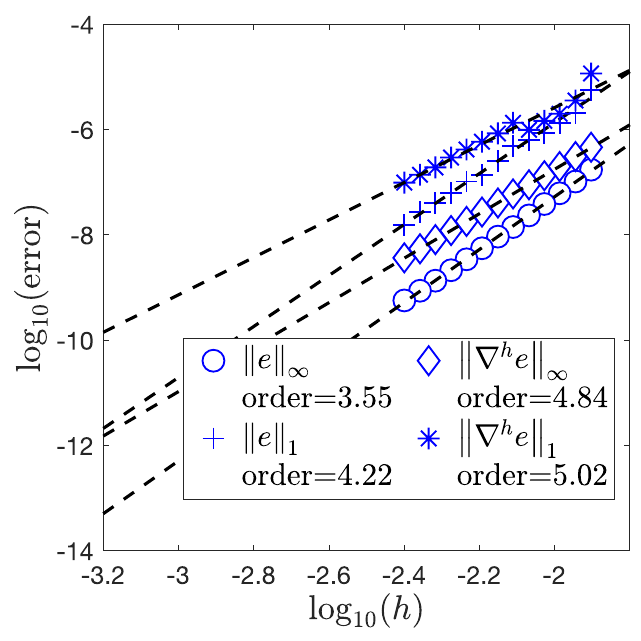}
\includegraphics[width=0.3\textwidth]{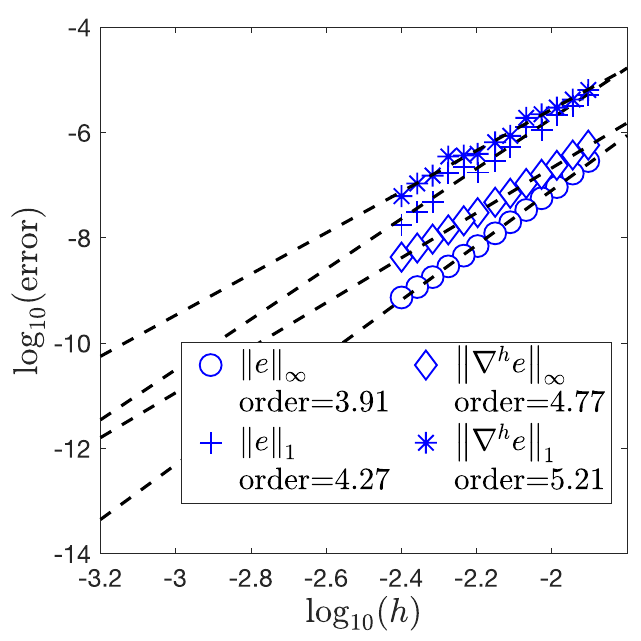}
\end{center}
\caption{Numerical errors for the test described in Section~\ref{sect:otherDomains}: the rotated leaf (left), the flower (middle), and the hourglass (right).}
	\label{fig:bestfit_domains}
\end{figure}

We display in Fig.~\ref{fig:bestfit_domains} the $L^1$ and $L^\infty$ errors along with the convergence curves.
The error decays monotonically and aligns closely with the optimal fourth-order reference slope, reaffirming  that the method is robust across different geometries and under mixed boundary conditions.

\subsection{Steady convection–diffusion problem with source term in an annulus}\label{sect:TestConvDiff}

Once again, the domain $\Omega$ is the annulus defined in Section~\ref{sect:test_hom}, and we seek $\phi$ solution of the linear convection-diffusion equation $-\kappa \Delta \phi+U\cdot \nabla \phi=g$, with $\kappa \neq 0$. We prescribe homogeneous Dirichlet boundary conditions to assess the ability of the scheme to catch boundary layers. To this end, the velocity has the following expression:
$$
U(x,y)=\frac{u_0}{r^2}\left (\begin{array}{c}x\\y\end{array}\right ),
$$
while $g(x,y)=g(r)=\frac{1}{r}$ depends only on $r$.
Therefore, the solution $\phi$ is invariant by rotation and the equation reads in polar coordinates
$$
-\frac{\kappa}{r}\frac{d}{dr}\left ( r\frac{d\phi}{dr}\right )+\frac{u_0}{r}\frac{d \phi}{dr}=\frac{1}{r}.
$$
Integration in $r$ yields $\displaystyle -\kappa r \frac{d\phi}{dr}+u_0\phi(r)=r+r_0,$ with $r_0\in \mathbb R$, and the homogeneous solution reads
$$
\tilde{\phi}_0(r) = A \phi_0(r)=Ar^{u_0/\kappa},\quad A\in\mathbb R.
$$
On the other hand, we assume that the general solution can be written as $\phi(r)=A(r) \phi_0(r)$ and get 
$$
-\kappa r^{1+u_0/\kappa}\frac{dA}{dr}=r+r_0.
$$
We have to consider several cases regarding the values of $u_0$ and $\kappa$.

\begin{description}
\item[Case 1: $u_0 \neq 0$ and $u_0 \neq \kappa$.] Integrating once again gets the relation
$$
A(r)=\frac{r^{1-u_0/\kappa}}{u_0-\kappa} +r_0\frac{r^{-u_0/\kappa}}{u_0}+C,
$$
and the general solution is given by
\begin{equation}\label{eq:exact_case1}
\phi(r)=A(r)\phi_0(r)=\frac{r}{u_0-\kappa}+\frac{r_0}{u_0}+Cr^{u_0/\kappa},
\end{equation}
where $C,\,r_0\in\mathbb R$ will be defined to satisfy the boundary condition. Assuming  $\phi(R_1)=\phi(R_2)=0$, we obtain the real constant values
$$
r_0=\frac{u_0}{\kappa-u_0}\ \frac{R_1 R_2^{u_0/\kappa}-R_2 R_1^{u_0/\kappa}}{R_2^{u_0/\kappa}-R_1^{u_0/\kappa}},\quad
C=\frac{1}{k-u_0}\frac{R_2-R_1}{R_2^{u_0/\kappa}-R_1^{u_0/\kappa}}.
$$

\item[Case 2: $u_0 = 0$.] For the pure diffusion problem, we get
$$
\phi(r)=A(r)\phi_0(r) = 
A(r)=-\frac{1}{\kappa} \left( r + r_0 \log \left| r \right| \right) + C, 
$$
where $r_0$ and $C$ are defined to satisfy the boundary condition. We assume that $\phi(R_1)=\phi(R_2)=0$ and deduce that
$$
r_0=-\frac{R_2-R_1}{\log(R_2)-\log(R_1)},\quad 
C=-\frac{R_2\log(R_1)-R_1\log(R_2)}{\kappa(\log(R_2)-\log(R_1))}.
$$
\item[Case 3: $u_0 = \kappa \neq 0$.] For this very specific case, we get an alternative solution that reads
$$
A(r)=-\frac{1}{\kappa} \left( \log \left| r \right| - \frac{r_0}{r} \right) + C
$$
and the general solution is then given by
$$
\phi(r)=-\frac{1}{\kappa} \left( r \log \left| r \right| - r_0 \right) + C r
$$
where $r_0$ and $C$ will be defined to satisfy the boundary condition. From $\phi(R_1)=\phi(R_2)=0$, we deduce that
$$
r_0=-\frac{R_1 R_2 (\log(R_2)-\log(R_1))}{R_2-R_1}, \quad C=\frac{R_2 \log(R_2)-R_1 \log(R_1)}{\kappa(R_2-R_1)}.
$$
\end{description}

\begin{figure}[ht]
\begin{center}
\includegraphics[width=0.3\textwidth]{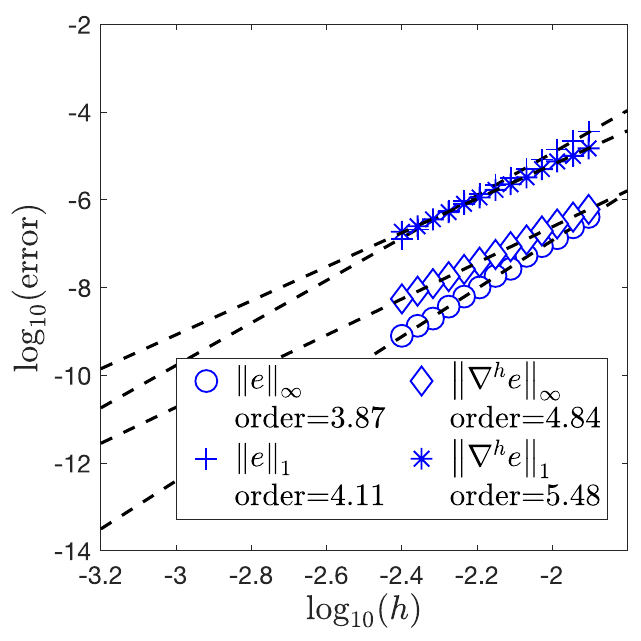}
\includegraphics[width=0.3\textwidth]{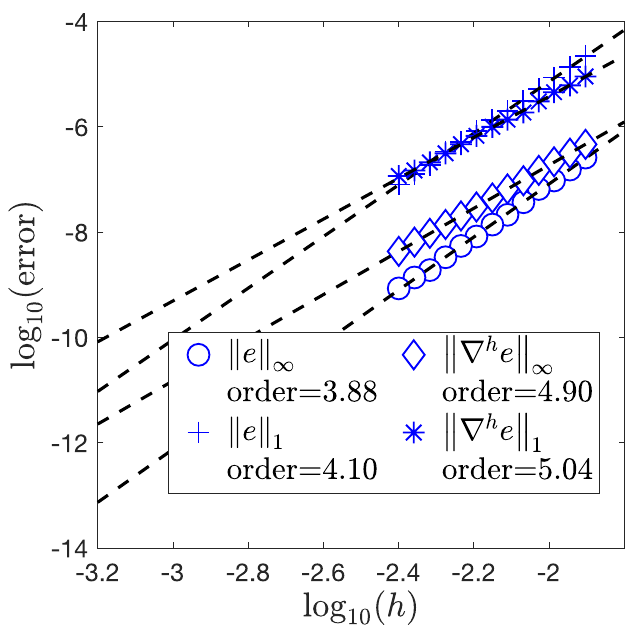}
\includegraphics[width=0.3\textwidth]{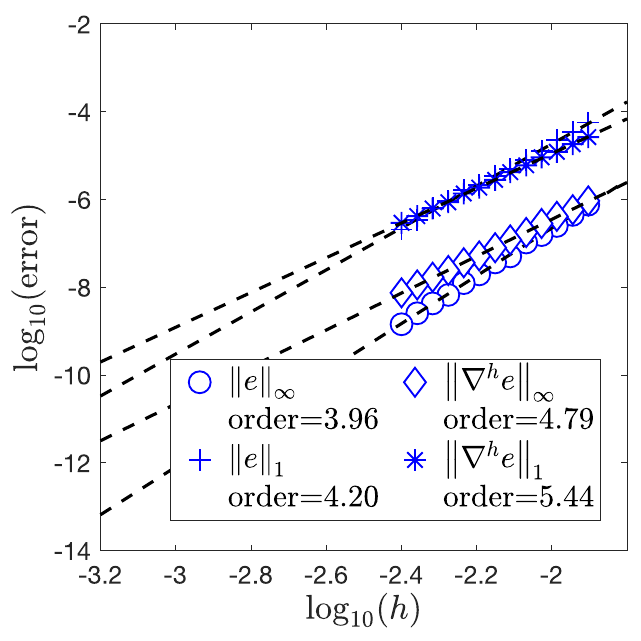}
\end{center}
\caption{
Numerical errors for the tests described in Section~\ref{sect:TestConvDiff}:
Case 1 with $k=2$ and $u_0=1$ (left),
Case 2 with $k=1$ and $u_0=0$ (middle), and
Case 3 with $k=1$ and $u_0=1$ (right).
}
\label{fig:TestConvDiff}
\end{figure}

We plot in Fig.~\ref{fig:TestConvDiff} the error convergence curves for the three cases. More specifically, the left panel corresponds to Case~1 with $\kappa = 2$ and $ u_0 = 1$, the middle panel displays Case~2 with $\kappa=1$ and $u_0=0$, and the right panel represents Case~3 with $\kappa=u_0=1$. The major conclusion is that we recover the optimal fourth order method in all the cases with very straightforward convergence curves and low deviation regarding the slope. This point is not so easy to achieve in the literature since the stencil configuration change for each grid and usually brings some instability/non monotony decay of the order, similar to the one we have observed with S1 or S2 in $L^\infty$ norm, for instance (Fig.~\ref{fig:BF_S123}).
\vskip 1em

\begin{figure}[ht]
\begin{center}
\includegraphics[width=0.3\textwidth]{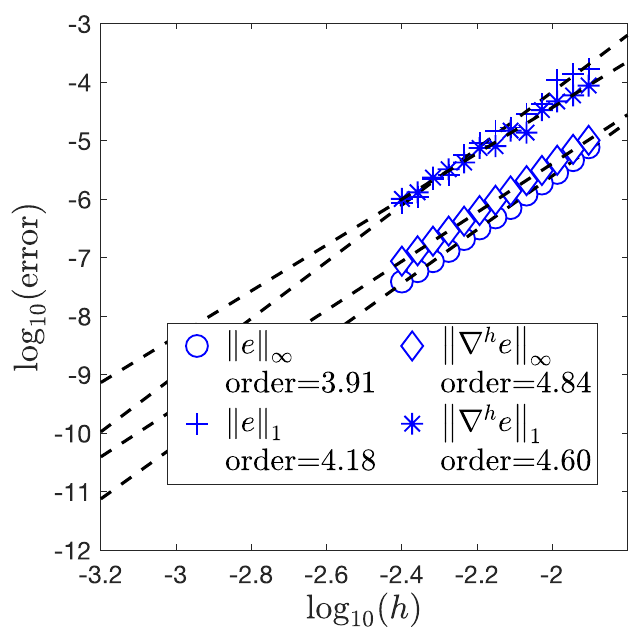}
\includegraphics[width=0.3\textwidth]{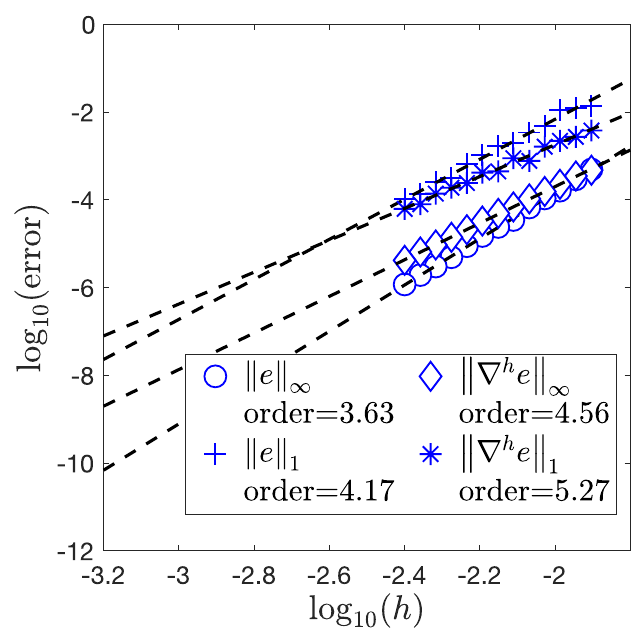}
\includegraphics[width=0.33\textwidth]{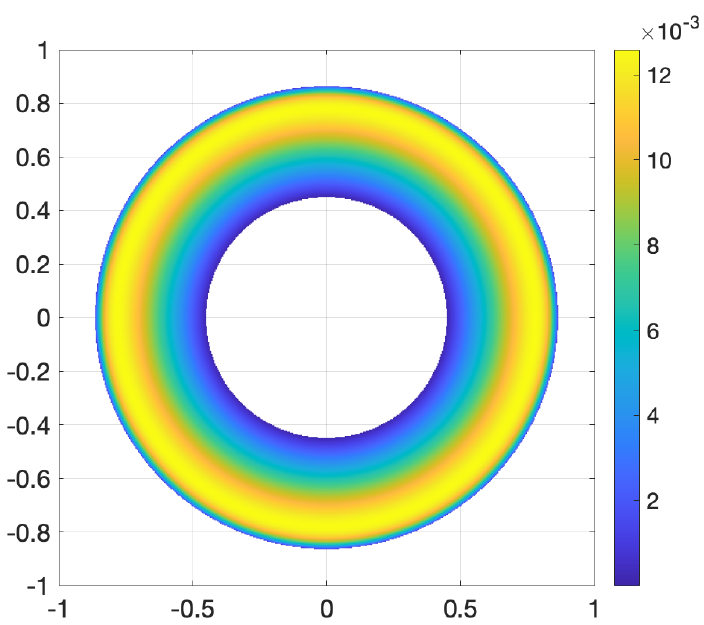}
\end{center}
\caption{
Numerical errors for the tests described in Section~\ref{sect:TestConvDiff}:
Benchmark 1: $\kappa=1$ and $u_0=10$ (left). Benchmark 2: $\kappa=1$ and $u_0=25$ (middle). The numerical solution of the case $\kappa=1$ and $u_0=25$ is shown in the right plot.
}
\label{fig:boundarylayer}
\end{figure}

Another issue concerns the Case~1, that exhibits pronounced boundary layer effects when \( u_0 R_2 \gg \kappa \neq 0 \) as shown in Fig.~\ref{fig:boundarylayer}, right panel. These arise from the dominance of the convective term over diffusion, combined with fixed Dirichlet boundary conditions, as evidenced by the analytical solution~\eqref{eq:exact_case1}.
The length $\delta=\frac{\kappa}{u_0}$ represents the characteristic size of the boundary layer to compare with the characteristic length of the geometry $R_2$, leading to the definition of the global dimensionless Péclet number $Pe=\frac{u_0 R_2}{\kappa}$. 
We perform two tests: benchmark 1 with $Pe=8$ and $\delta=1/10$, and benchmark 2 with $Pe=20$ and $\delta=1/25$.
Since our objective is to investigate how the boundary conditions' discretization performs in the presence of such boundary layers, 
we introduce the dimensionless local (cell) Péclet number, defined as $Pe_\text{loc}=\frac{u_0 \Delta x}{\kappa}$, which depends on the grid size. Such a number describes the trade-off between boundary layer and grid size at the discrete level.

Taking $\Delta x=\frac{2}{N}$ with $N$ ranging from 160 to 502, we have $Pe_\text{loc}\in [4.0\times 10^{-2},1.6\times10^{-1}]$ for the first benchmark (Fig.~\ref{fig:boundarylayer}, left panel), and $Pe_\text{loc}\in [1.0\times10^{-2},4.2\times10^{-1}]$ for the second benchmark (Fig.~\ref{fig:boundarylayer}, middle panel), where we report the numerical errors and convergence orders. Once again, we observe the errors decrease monotonically and in accordance with the expected fourth-order convergence rate. We point out that for problems with a higher cell Péclet number, namely in convection-dominated regimes, a stable and accurate discretization of the convective term would typically require an upwind scheme rather than central differences, unless a sufficiently fine grid is employed to satisfy the Péclet condition.

\section{Conclusive Remarks}\label{sect:concl}
In this work, we introduced a novel boundary operator framework for ghost-point finite-difference methods on unfitted Cartesian grids, with the primary goal of achieving high-order accuracy while maintaining compact and well-conditioned stencils. The proposed approach reformulates the enforcement of boundary conditions through a local linear operator constructed via constrained least-squares reconstruction. This strategy provides significant flexibility in stencil design. 

A key contribution of the paper lies in the identification of two criteria, which guide the construction of admissible and stable stencils: local and global conditioning. Building upon these criteria, we developed several strategies, finally leading to the cone-based adaptive procedure. This approach successfully balances the requirements of admissibility, conditioning, and compactness, overcoming limitations of traditional wide-stencil ghost-point methods.

The numerical experiments confirm several important findings. First, the proposed method achieves consistent fourth-order accuracy for both the solution and its gradient across a wide range of test cases, including mixed boundary conditions and complex geometries represented implicitly via level-set functions. Second, unlike many existing ghost-point techniques, the method exhibits smooth and monotonic convergence behaviour with respect to grid refinement, reducing the occurrence of error oscillations typically associated with irregular stencil configurations. Third, the compactness of the selected stencils makes the method particularly suitable for large-scale and parallel simulations, where communication overhead is a critical factor.

Moreover, the method demonstrates robustness in challenging regimes, including convection-diffusion problems with boundary layers. Even in such cases, the approach maintains stability and accuracy without requiring excessive stencil enlargement. 

Overall, the boundary operator formulation provides a flexible framework for high-order discretizations on unfitted domains, contributing to the design of efficient and robust numerical methods on complex geometries, particularly interesting in applications involving moving interfaces or multiscale features.

\section*{Acknowledgments}

\noindent This work was partially supported by: 
\begin{itemize}
\item MUR Prin 2022, project code 2022N9BM3N, Grant No. CUP: E53D23005830006;
\item MUR Prin PNRR 2022, project code P2022BNB97, Grant No. CUP: E53D23017960001;
\item Spoke 1 Future HPC \& Big Data of the Italian Research Center on High-Performance Computing, Big Data and Quantum Computing (ICSC) funded by MUR Missione 4 Componente 2 Investimento 1.4: Potenziamento strutture di ricerca e creazione di “campioni nazionali di R \&S (M4C2-19)” - Next Generation EU (NGEU);
\item MUR Prin 2022, project code 2022M9BKBC, Grant No. CUP  D53D23005880006;
\item MUR Prin PNRR 2022, project code P202254HT8, Grant No. CUP B53D23027760001.
\end{itemize}

\noindent Armando Coco and Alessandro Coclite are members of Gruppo Nazionale per il Calcolo Scientifico (GNCS) of the Istituto Nazionale di Alta Matematica (INdAM).\\

\noindent The authors Stéphane Clain and Rui Miguel Pereira were financed by national funds from FCT - Fundação para a Ciência e a Tecnologia, I.P., within the context of the project HiPMos4C setting with reference 2023.16854.ICDT2023.16854.ICDT, https://doi.org/10.54499/2023.16854.ICDT.\\

\noindent Stéphane Clain would like to acknowledge the support of Centro de Matemática of the Coimbra University (CMUC), funded by the Portuguese Foundation for Science and Technology, under the references  UIDB/00324/2020 (https://doi.org/10.54499/UIDB/00324/2020) and UIDP/00324/2020 (https://doi.org/10.54499/UIDP/00324/2020)\\

\noindent Rui Pereira was supported by national funds through FCT  - Fundação para a Ciência e a Tecnologia, under the Program Contract of the Centre of Mathematics of the University of Minho (CMAT/UM), UID/00013/2025 (https://doi.org/10.54499/UID/00013/2025).

\section*{Conflicts of Interest} The authors declare that there is no conflict of interest.

\bibliographystyle{unsrtnat}
\bibliography{biblio.bib}

\appendix
\section{The S4 algorithm.}\label{sec::appA}
\begin{algorithm}[H]
\caption{Construction of the S4.1, S4.2 and S4.3 isosceles triangular stencils (vertex angle $\vartheta$)}
\label{alg:S4}
\begin{algorithmic}[1]
\State Define the set of candidate grid points $\by_k$
such that the angle between 
$\overline{\bx_k\by_k}$ and $\overline{\bx_k\bp_k}$ are $\le \vartheta/2$ and order them by the closest to the ghost point:
$\{\by_1, \by_2, \by_3, \ldots \}$ 
\State Define $q=N_o$ and the stencil $\stencil^k = \{\by_1, \ldots, \by_{q} \}$
\While{$\rho(B_k)>10^6$}
    \State $\stencil^k \leftarrow \stencil^k \cup \{ \by_{q+1} \}$
    \State $q \leftarrow q+1$
\EndWhile
\State Set $\stencil^k_\text{S4.1} = \stencil^k$
\If{$\rho(B_k)>10^6$}
    \State Identify the point $\bx_\ell \in \stencil^k$ corresponding to the largest coefficient $\ba^\ell$ (Def.~\ref{defn:wellcond}).
    \State $\stencil^k \leftarrow \stencil^k \backslash \{ \bx_\ell \}$.
    \State $\stencil^k \leftarrow \stencil^k \cup \{ \by_{q+1} \}$
\EndIf
\State Set $\stencil^k_\text{S4.2} = \stencil^k$
\If{$\rho(B_k)>10^6$}
    \State Modify the collar point $\bp_k$ by projecting the ghost point horizontally or vertically onto the boundary (choose the closer projection if both exist).
    \State Rebuild stencil $\stencil^k$ using Strategy~2.
\EndIf
\State Set $\stencil^k_\text{S4.3} = \stencil^k$
\end{algorithmic}
\end{algorithm}
\end{document}